\documentclass[]{article}
\usepackage{graphicx}
\usepackage{amsfonts}
\usepackage{amsmath}
\usepackage{amssymb}
\usepackage{fancyhdr}
\usepackage{titlesec}
\usepackage{indentfirst}
\usepackage{booktabs}
\usepackage{verbatim}
\usepackage{color}
\usepackage{amsthm}
\usepackage{hyperref} 
\usepackage{subcaption}

\usepackage[page,header]{appendix}
\usepackage{titletoc}

\numberwithin{equation}{section}
\newtheorem{theorem}{Theorem}[section]

\newtheorem{remark}[theorem]{Remark}

\newcommand{\rd}{\,\mathrm{d}}
\newcommand{\mM}{\mathcal{M}}

\makeatletter
\usepackage{fullpage}
\usepackage{url}
\title{Asymptotic-Preserving Dynamical Low-Rank Method \\for the Stiff Nonlinear Boltzmann Equation}
\author{Lukas Einkemmer\footnote{Department of Mathematics, University of Innsbruck, Innsbruck, A-6020, Austria (lukas.einkemmer@uibk.ac.at).}, \ Jingwei Hu\footnote{Department of Applied Mathematics, University of Washington, Seattle, WA 98195, USA (hujw@uw.edu).}, 
	\ and \  Shiheng Zhang\footnote{Department of Applied Mathematics, University of Washington, Seattle, WA 98195, USA (shzhang3@uw.edu).} }  
\makeatother
\usepackage{babel}
\begin{document}
\maketitle
\begin{abstract}
In kinetic theory, numerically solving the full Boltzmann equation is extremely expensive. This is because the Boltzmann collision operator involves a high-dimensional, nonlinear integral that must be evaluated at each spatial grid point and every time step. The challenge becomes even more pronounced in the fluid (strong collisionality) regime, where the collision operator exhibits strong stiffness, causing explicit time integrators to impose severe stability restrictions. In this paper, we propose addressing this problem through a dynamical low-rank (DLR) approximation. The resulting algorithm requires evaluating the Boltzmann collision operator only $r^2$ times, where $r$, the rank of the approximation, is much smaller than the number of spatial grid points. We propose a novel DLR integrator, called the XL integrator, which reduces the number of steps compared to the available alternatives (such as the projector splitting or basis update \& Galerkin (BUG) integrator). For a class of problems including the Boltzmann collision operator which enjoys a separation property between physical and velocity space, we further propose a specialized version of the XL integrator, called the sXL integrator. This version requires solving only one differential equation to update the low-rank factors. Furthermore, the proposed low-rank schemes are asymptotic-preserving, meaning they can capture the asymptotic fluid limit in the case of strong collisionality. Our numerical experiments demonstrate the efficiency and accuracy of the proposed methods across a wide range of regimes, from non-stiff (kinetic) to stiff (fluid).

\end{abstract}

{\small 

}

\section{Introduction}
\label{sec:intro}

The Boltzmann equation \cite{Cercignani} is one of the fundamental equations in kinetic theory. In non-dimensionalized form, the equation reads as
\begin{equation} \label{boltz}
     \partial_t f +v \cdot \nabla_x f =\frac{1}{\varepsilon}\mathcal{Q}(f, f), \quad t>0, \  x\in \Omega_x\in \mathbb{R}^{d_x}, \  v\in \mathbb{R}^{d_v},
 \end{equation}
where $f(x,v,t)$ is the phase-space distribution function of position $x$, velocity $v$, and time $t$. $\mathcal{Q}$ is the Boltzmann collision operator modeling binary collisions between particles, whose bilinear form is given by
\begin{equation}
\mathcal{Q}(g,f)(v)=\int_{\mathbb{R}^{d_v}}\int_{S^{d_v-1}} B(|v-v_*|,\cos \theta)[g(v_*')f(v')-g(v_*)f(v)]\rd{\sigma}\rd{v_*}, \label{collision-integral}
\end{equation}
where the post-collisional velocities $(v',v_*')$ are related to pre-collisional velocities $(v,v_*)$ through the conservation of momentum and energy:
\begin{equation}
v'=\frac{v+v_*}{2}+\frac{|v-v_*|}{2}\sigma, \quad v_*'=\frac{v+v_*}{2}-\frac{|v-v_*|}{2}\sigma,
\end{equation}
with $\sigma$ being a vector on the unit sphere $S^{d_v-1}$. The collision kernel $B$ is a non-negative function encoding the collision frequency, which depends on the relative velocity $|v-v_*|$ and cosine of the deviation angle $\theta$, i.e., the angle between $v-v_*$ and $v'-v_*'$. Finally, $\varepsilon$ is the Knudsen number, defined as the ratio of the mean free path to the characteristic length. The value of $\varepsilon$ indicates the regime of the problem: for $\varepsilon=\mathcal{O}(1)$, the system is in the kinetic regime, whereas for $\varepsilon\ll 1$, the system approaches the fluid regime. In this latter case, one can derive the fluid limits of the Boltzmann equation \cite{BGL91}.

The collision operator satisfies the following conservation property:
\begin{equation}
\int_{\mathbb{R}^{d_v}} \mathcal{Q}(f,f)\phi(v)\rd{v}=0, \quad \phi(v):=(1,v,|v|^2/2)^T.
\end{equation}
Furthermore, the kernel of the collision operator is given by the Maxwellian $\mathcal{M}[f]$, 
\begin{equation} \label{Max}
\mathcal{Q}(f,f)=0 \ \Longleftrightarrow \ f=\mathcal{M}[f]:=\frac{\rho(x,t)}{(2\pi T(x,t))^{d_v/2}} e^{-\frac{|v-u(x,t)|^2}{2T(x,t)}},
\end{equation}
where \( \rho \), \( u \), and \( T \) are the density, bulk velocity, and temperature, defined through the moments of $f$:
\begin{equation} \label{macro-quantity}
\rho = \int_{\mathbb{R}^{d_v}} f\rd{v}, \quad 
u = \frac{1}{\rho} \int_{\mathbb{R}^{d_v}} f v \rd{v}, \quad 
T = \frac{1}{d_v\rho} \int_{\mathbb{R}^{d_v}} f |v-u|^2 \rd{v}.
\end{equation}
It is easy to see that $f$ and $\mathcal{M}[f]$ share the first $d_v+2$ moments, i.e.,
\begin{equation}
\int_{\mathbb{R}^{d_v}}f\phi(v)\rd{v}=\int_{\mathbb{R}^{d_v}}\mathcal{M}[f]\phi(v)\rd{v}=(\rho,\rho u,E)^T:=\mathcal{U}, 
\end{equation}
where $E:=\frac{d_v}{2}\rho T+\frac{1}{2}\rho u^2$ is the total energy.

Using the above properties, one can formally derive the compressible Euler equations as the leading-order limit of the Boltzmann equation \eqref{boltz} as $\varepsilon \rightarrow 0$: taking the moments $\int_{\mathbb{R}^{d_v}} \cdot \,\phi(v) \rd{v}$ on both sides of \eqref{boltz} yields the local conservation laws:
\begin{equation} \label{local-con}
\partial_t \mathcal{U}+\nabla_x \cdot \int_{\mathbb{R}^{d_v}}f v\phi(v)\rd{v}=0.
\end{equation}
The system \eqref{local-con} is not closed. However, when $\varepsilon \rightarrow 0$, \eqref{boltz} implies that $\mathcal{Q}(f,f)\rightarrow 0$, hence $f\rightarrow \mathcal{M}[f]$. Then \eqref{local-con} can be closed (the flux term can be completely represented by $\rho$, $u$, and $T$), yielding the compressible Euler equations:
\begin{align}  \label{Euler}
\left\{
\begin{array}{l}
\displaystyle \partial_t \rho+\nabla_x\cdot (\rho u)=0,\\[8pt]
\displaystyle \partial_t (\rho u)+\nabla_x\cdot (\rho u\otimes u) +\nabla_x(\rho T)=0,\\[8pt]
\displaystyle \partial_t E+\nabla_x\cdot ((E+\rho T) u)=0.
\end{array}\right.
\end{align}

Solving the compressible Euler equations numerically is orders of magnitude cheaper than solving the Boltzmann equation. This is primarily due to the fact that the phase space $(x,v)$ can be up to six-dimensional,  and for each point in phase space, one must compute \eqref{collision-integral}, a $(2d_v-1)$ dimensional collisional integral. Thus, a naive algorithm would scale as $\mathcal{O}(N^{d_x}MN^{2d_v})$, where $N$ is the number of grid points in each direction of the phase space, and $M$ is the number of points used on the sphere $S^{d_v-1}$. 

Much work in the literature has focused on improving the collision solver. For example, a direct Fourier spectral method was introduced in \cite{PR00}, reducing the complexity of a one-time evaluation of the collision integral from $\mathcal{O}(MN^{2d_v})$ to $\mathcal{O}(N^{2d_v})$. Later, fast Fourier spectral methods were introduced in \cite{MP06} for certain collision kernels, achieving a complexity of $\mathcal{O}(MN^{d_v}\log N)$, and in \cite{GHHH17} for general collision kernels, achieving a complexity of $\mathcal{O}(MN^{d_v+1}\log N)$ (note that in both of these works, $M$ can be chosen much smaller than $N^{d_v-1}$). Despite these advances, numerically solving the full spatially inhomogeneous Boltzmann equation remains computationally prohibitive because the above costs must be multiplied by $N^{d_x}$.

Thus, our goal in this paper is to remove the multiplicative dependence on $N^{d_x}$. More specifically, we aim to reduce the number of dimensions (and thus the memory complexity) required to store the phase-space distribution $f$. An added benefit of this is a drastic reduction in the number of calls to the fast collision solver, which in a classic algorithm must be evaluated for every spatial grid point. To accomplish this, we will use a dynamical low-rank (DLR) method. 

DLR approximations trace their history back to the famous method of Hartree and Fock, or to be more precise the multi-configuration time dependent Hartree (MCTDH) method  (see, e.g., \cite{Meyer2009,lubich2008}). The main idea is that the phase-space distribution function is approximated by a low-rank representation
\[ f(x,v,t) = \sum_{i=1}^r X_i(x,t) S_{ij}(t) V_j(v,t),\]
where $r$ is the rank of the approximation. For the low-rank factors $X_i$, $V_j$, and $S_{ij}$, evolution equations are then derived. It has been shown that for many collisionless (see, e.g., \cite{Einkemmer2018,Einkemmer2020,Einkemmer2023,Guo2022b}) and collisional (see, e.g., 
\cite{peng2020,Kusch2021a, HW22, CHS24}) kinetic problems a relatively small rank is required to accurately represent the phase-space distribution function. This then drastically reduces the storage cost from $\mathcal{O}(N^{d_x+d_v})$ to $\mathcal{O}(r(N^{d_x} + N^{d_v}))$. Such low-rank methods have been previously considered for simplified collision operators such as models related to linear radiative transfer \cite{peng2020,Ding2021,Einkemmer2021,Baumann2023}, Bhatnagar–Gross–Krook (BGK) \cite{Einkemmer2021d,Dektor2024}, and Fokker–Planck type operators \cite{Coughlin2022,Nakao2023}. In these situations we are often interested in numerical methods that can capture the asymptotic limit. In the context of the Boltzmann equation, this means a scheme that solves the equation \eqref{boltz}, and automatically reduces to a consistent discretization to the compressible Euler limit \eqref{Euler} as $\varepsilon \rightarrow 0$, i.e., asymptotic-preserving (AP) \cite{Jin_Rev, HJL17}. Usually, being AP requires some sort of implicit treatment of at least a part of the equations. Although, for many of the simplified collision operators mentioned above schemes can be devised in such a way that the implicit update can be computed very efficiently, for example, the IMEX schemes in \cite{Einkemmer2021,Einkemmer2022a} for the linear radiative transfer equation and \cite{Einkemmer2021d} for the BGK equation.

The purpose of this paper is to derive a dynamical low-rank method for the stiff Boltzmann equation \eqref{boltz} that preserves the compressible Euler limit \eqref{Euler} as $\varepsilon\rightarrow 0$ (i.e.,~is AP). This adds a number of challenges. In particular, being able to use a fast collision solver is essential for computational efficiency and we want to reduce the number of times this solver is called as much as possible. The proposed algorithm requires no modification of the fast collision solver, and calls it only $\mathcal{O}(r^2) \ll \mathcal{O}(N^{d_x})$ many times. To treat the stiffness and obtain an AP scheme we propose a splitting between advection and collision operators and then penalize the Boltzmann collision operator. The penalization is treated implicitly (within an IMEX scheme), but as it is much simpler than the Boltzmann collision operator this can be done easily within the DLR framework. The splitting between advection and collision has the further advantage that the DLR schemes we propose (called the XL integrator and sXL integrator) are much simpler and thus more efficient than what would be obtained by applying one of the standard robust integrators (such as the projector splitting integrator \cite{lubich2014projector} or the BUG integrator \cite{ceruti2022rank}).

The rest of this paper is organized as follows. Although the XL and sXL integrators are designed with the Boltzmann equation in mind, they are in fact quite generic and can be applied to general matrix differential equations. For this reason, we start in Section~\ref{sec:general XL integrator} to describe these integrators in a general setting. We then introduce the proposed DLR methods for the Boltzmann equation, first in the non-stiff case (Section~\ref{sec:nonstiff-Boltz}) and then in the stiff case (Section~\ref{sec:stiff-Boltz}). In Section~\ref{sec:BGK operator and proof of AP property}, we prove the AP property of the proposed methods for a special BGK equation. Numerical examples are presented in Section~\ref{sec:experiments} to demonstrate the efficiency and accuracy of the proposed methods across a wide range of regimes, from non-stiff (kinetic) to stiff (fluid). The paper is concluded in Section~\ref{sec:conclusions}.


\section{The general XL and sXL integrators}
\label{sec:general XL integrator}

In the dynamical low-rank literature a number of robust integrators have been proposed to solve partial differential equations. Among them are the projector-splitting or KSL integrator \cite{lubich2014projector}, the family of BUG integrators \cite{ceruti2022unconventional,ceruti2022rank}, and the parallel integrator \cite{Ceruti2023}. For more information on different robust integrators, we refer to the recent review article on low-rank methods \cite{Einkemmer2024c}. Commonly, these integrators are introduced for spatially discretized and time-continuous matrix equations, but they can also be formulated in the fully continuous case \cite{Einkemmer2018}. We now introduce the new XL integrator and sXL integrator for a general matrix differential equation, before we turn our attention to the Boltzmann equation in the next two sections.

Consider a $m_1\times m_2$ matrix differential equation
\begin{equation} \label{matrix-eq}
\dot{A}(t)=F(t,A(t)). 
\end{equation}
We start from a factored rank-$r$ matrix $Y^n=X^nS^n(V^n)^T$ at time $t^n$ and propose an {\bf XL integrator} to obtain the solution $Y^{n+1}=X^{n+1}S^{n+1}(V^{n+1})^T$ at time $t^{n+1}=t^n+\Delta t$ with rank $r$ as follows: 
\begin{itemize}
\item \textbf{X-step}: Integrate from $t^n$ to $t^{n+1}$ the $m_1\times r$ matrix equation
\begin{equation}
\dot{K}(t)=F(t,K(t)(V^n)^T)V^n, \quad K(t^n)=X^nS^n,
\end{equation}
to obtain $K(t^{n+1})$. Determine $\hat{X}\in \mathbb{R}^{m_1\times 2r}$ as
\begin{equation}
\hat{X}=\text{orth}([X^n,K(t^{n+1})]),
\end{equation}
where $\text{orth}(\cdot)$ denotes an orthonormalization algorithm (e.g., a QR decomposition computed by Householder reflections or using the modified Gram-Schmidt process). We assume here and in the following that since $X^n$ is already orthogonal it is left invariant by the orthogonalization, that is, $\hat{X} = [X^n, \dots]$.
\item \textbf{L-step}: Integrate from $t^n$ to $t^{n+1}$ the $m_2\times 2r$ matrix equation
\begin{equation}
\dot{\hat{L}}(t)=F(t,\hat{X}\hat{L}(t)^T)^T\hat{X}, \quad \hat{L}(t^n)=[V^n(S^n)^T,0],
\end{equation}
to obtain $\hat{L}(t^{n+1})$, where $0$ is the $m_2\times r$ zero matrix. Determine $\hat{V}\in \mathbb{R}^{m_2\times 2r}$, $\hat{S}\in\mathbb{R}^{2r\times 2r}$ as
\begin{equation}
[\hat{V},\hat{S}]=QR(\hat{L}(t^{n+1})).
\end{equation}
This way, we obtain the solution $\hat{Y}^{n+1}=\hat{X}\hat{S}\hat{V}^T$ with rank $2r$.
\item \textbf{Truncation step}: Use the singular value decomposition (SVD) of $\hat{S}$ to reduce the rank back to $r$ and to obtain the solution $Y^{n+1}=X^{n+1}S^{n+1}(V^{n+1})^T$ at the next time step.
\end{itemize}

\begin{remark}
Note that the above truncation step is exactly the same as that in the augmented BUG integrator \cite{ceruti2022rank}. Just like the rank adaptivity in that approach, it is not necessary to truncate the rank back to $r$; instead, one can choose the rank adaptively by setting a certain threshold. For simplicity, we adhere to the former choice in this paper.
\end{remark}

Like the above-mentioned existing DLR integrators, the proposed \textbf{XL integrator} enjoys the exactness property.
\begin{theorem}[\textbf{Exactness property}]
        Let \( A(t) \in \mathbb{R}^{m_1 \times m_2} \) have rank \( r \) for all \( t \in [t^n, t^{n+1}] \), then it admits the factorization
    \(
        A(t) = X(t) L(t)^{T} = X(t) S(t) V(t)^{T}.
    \)
    Assume that \( V(t^{n+1})^{T} V(t^{n}) \) is invertible. With \( Y^{n} = A(t^{n}) \), the \textbf{XL integrator} for \eqref{matrix-eq} is exact, that is, $Y^{n+1} = A(t^{n+1})$, provided that at least \( r \) singular values are retained during the truncation step.
\end{theorem}
\begin{proof}
    In the {\bf X-step}, we can choose w.l.o.g. that $\hat{X}=[X^n,\tilde{X}^{n+1}]$, where $\tilde{X}^{n+1}$ is the orthogonal factor from $K(t^{n+1})$ and $X^n(\tilde{X}^{n+1})^T = 0$. \cite[Lemma 1]{ceruti2022unconventional} shows that $\tilde{X}^{n+1}$ and $A(t^{n+1})$ have the same range, or equivalently,
    \begin{equation*}
        \tilde{X}^{n+1}(\tilde{X}^{n+1})^TA(t^{n+1}) = A(t^{n+1}).
    \end{equation*}
    We then also have 
    \begin{equation*}
        \hat{X}(\hat{X})^TA(t^{n+1}) = A(t^{n+1}).
    \end{equation*}
    Hence,
    \begin{equation*}
    \begin{aligned}
         \hat{X}\hat{S}\hat{V}^T &= \hat{X}\hat{L}^T\\ &= \hat{X}\left(\hat{L}(t^n)^T + (\hat{X})^T(A(t^{n+1})-A(t^n)) \right)\\
         &= \hat{X}(\hat{X})^TA(t^{n+1}) + 
         \hat{X}\hat{L}(t^n)^T -\hat{X}(\hat{X})^TA(t^n)\\
         &=A(t^{n+1}) +{X^n}{L}(t^n)^T -\hat{X}(\hat{X})^TA(t^n)\\
         &=A(t^{n+1}).
    \end{aligned}
    \end{equation*}
   Given that the matrix is already of rank \( r \), a truncation to rank \( r \) has no effect, ensuring the integrator is exact.
\end{proof}

Let us then assume the right hand side of the equation \eqref{matrix-eq} is autonomous and $F(A(t))$ satisfies the following {\it separable} property: if
\begin{equation}
Y=KV^T, \quad K\in\mathbb{R}^{m_1\times r}, \quad V\in \mathbb{R}^{m_2\times r}, 
\end{equation}
then
\begin{equation}
F(Y) =F(KV^T)= \mathcal{K}(K)\mathcal{V}(V)^T, \quad \mathcal{K}: \mathbb{R}^{m_1\times r} \rightarrow \mathbb{R}^{m_1\times \tilde{r}}, \ \mathcal{V}: \mathbb{R}^{m_2\times r} \rightarrow \mathbb{R}^{m_2\times \tilde{r}},
\end{equation}
where typically $\tilde{r}\geq r$. With this property, the {\bf X-step} of the above {\bf XL integrator} becomes:
integrate from $t^n$ to $t^{n+1}$ the $m_1\times r$ matrix equation
\begin{equation}
\dot{K}(t)=F(K(t)(V^n)^T)V^n=\mathcal{K}(K(t))\mathcal{V}(V^n)^TV^n, \quad K(t^n)=X^nS^n,
\end{equation}
to obtain $K(t^{n+1})$. Then it is clear that the effect of this step, if a forward Euler scheme is applied, is simply to add $\mathcal{K}(K(t^{n}))$ to the basis $X^n$. Thus, instead of solving the equation for $K$, we propose the following first-order {\bf sXL integrator} (XL integrator for separable problems):
\begin{itemize}
\item \textbf{Augmentation step}: Determine $\hat{X}\in \mathbb{R}^{m_1\times (r+\tilde{r})}$ as
\begin{equation}
\hat{X}=\text{orth}([X^n,\mathcal{K}(X^nS^n)]).
\end{equation}
\item \textbf{L-step}: Solve the $m_2\times (r+\tilde{r})$ matrix equation
\begin{equation}
\frac{\hat{L}^{n+1}-\hat{L}^n}{\Delta t}=F(t,\hat{X}(\hat{L}^n)^T)^T\hat{X}, \quad \hat{L}^n=[V^n(S^n)^T,0],
\end{equation}
where $0$ is the $m_2\times \tilde{r}$ zero matrix. Determine $\hat{V}\in \mathbb{R}^{m_2\times (r+\tilde{r})}$, $\hat{S}\in\mathbb{R}^{(r+\tilde{r})\times (r+\tilde{r})}$ as
\begin{equation}
[\hat{V},\hat{S}]=QR(\hat{L}^{n+1}).
\end{equation}
This way, we obtain the solution $\hat{Y}^{n+1}=\hat{X}\hat{S}\hat{V}^T$ with rank $r+\tilde{r}$.

\item \textbf{Truncation-step}: Use the SVD of $\hat{S}$ to reduce the rank back to $r$ and to obtain $Y^{n+1}=X^{n+1}S^{n+1}(V^{n+1})^T$.
\end{itemize}

\begin{remark}
Similar to the {\bf XL integrator} and {\bf sXL integrator}, one can do an analogous {\bf VK integrator} and {\bf sVK integrator} to simplify the V-step. Which one to use depends on the structure of the specific equation.
\end{remark}


\section{Dynamical low-rank method for the non-stiff Boltzmann equation}\label{sec:Dynamical low rank methods for the Boltzmann equation}
\label{sec:nonstiff-Boltz}

In this section, we briefly review the dynamical low-rank method for the non-stiff Boltzmann equation (i.e., equation \eqref{boltz} with $\varepsilon=1$) as proposed in \cite{hu2022adaptive}. The method is a straightforward application of the projector-splitting integrator \cite{lubich2014projector}. However, an important feature we wish to highlight is that the collision operator $\mathcal{Q}(f,f)$, despite its complicated structure, exhibits some favorable properties within the low-rank framework, which will also be important for our later discussion in the stiff case:
\begin{itemize}
\item $\mathcal{Q}$ is a bilinear operator, whose nonlinearity is benign compared to some simplified collision models, such as the BGK operator (see Section~\ref{sec:stiff-Boltz} for more discussion); 
\item $\mathcal{Q}$ is an operator acting only in the velocity space, with the spatial variable serving merely as a parameter. 
\end{itemize}

We first write the equation \eqref{boltz} with $\varepsilon=1$ as
\begin{equation}\label{Boltzmann-rhs}
     \partial_t f = -v \cdot \nabla_x f + \mathcal{Q}(f, f) := \text{RHS}.
 \end{equation}
We then constrain the distribution function $f(x,v,t)$ to a low-rank manifold $\mathcal{M}_r$ such that
\begin{equation}\label{fsvd}
    f(x, v, t) = \sum_{i,j=1}^r X_i(x, t) S_{ij}(t) V_j(v, t),
\end{equation}  
where \( r \) represents the rank of the approximation, and the basis functions \( \{X_i\}_{i=1}^r \subset L^2(\Omega_x) \) and \( \{V_j\}_{j=1}^r \subset L^2(\Omega_v) \) are orthonormal: 
\begin{equation}\label{orth}
    \langle X_i, X_j \rangle_x = \delta_{ij}, \quad \langle V_i, V_j \rangle_v = \delta_{ij}, \quad 1 \leq i, j \leq r,
\end{equation}  
with \( \langle \cdot, \cdot \rangle_x \) and \( \langle \cdot, \cdot \rangle_v \) denoting the inner products in \( L^2(\Omega_x) \) and \( L^2(\Omega_v) \), respectively. Here \( \Omega_v \subset \mathbb{R}^{d_v} \) is a truncated velocity domain, which is a reasonable assumption when performing numerical simulations of kinetic equations. 

By projecting \eqref{Boltzmann-rhs} onto the tangent space of $\mathcal{M}_r$, one obtains
\begin{equation}\label{projection}
     \partial_t f = P_f(\text{RHS}),
\end{equation}
where the orthogonal projector \( P_f \) can be expressed as
\begin{equation} 
    P_f(\text{RHS}) = \sum_{j=1}^r \langle V_j, \text{RHS} \rangle_v V_j 
- \sum_{i,j=1}^r X_i \langle X_i V_j, \text{RHS} \rangle_{x,v} V_j 
+ \sum_{i=1}^r X_i \langle X_i, \text{RHS} \rangle_x.
\end{equation}
To avoid the inversion of a potentially singular matrix arising from this formulation, the projector-splitting integrator splits the right-hand side of \eqref{projection} into three subprojections,
and solves each subproblem consecutively. Since the equation \eqref{Boltzmann-rhs} is non-stiff, applying a fully explicit first-order discretization in each subproblem results in the following {\bf KSL integrator}: 

Given $\{X_i^n\}_{i=1}^r$, $\{V_j^n\}_{j=1}^r$, $\{S_{ij}^n\}_{i,j=1}^r$ at time $t^n$, it produces $\{X_i^{n+1}\}_{i=1}^r$, $\{V_j^{n+1}\}_{j=1}^r$, $\{S_{ij}^{n+1}\}_{i,j=1}^r$ at $t^{n+1}=t^n+\Delta t$ as follows.
\begin{itemize} 
    \item 
    \textbf{K-step}: 
Define $K_j^{n} = \sum_{i=1}^{r} X_i^{n} S_{ij}^{n}$, and solve for $K_j^{(1)}$ using the equation
\begin{align} \label{K-scheme}
     \frac{K_{j}^{(1)}-K_{j}^{n}}{\Delta t}=-\sum_{l=1}^{r}\langle vV_{j}^{n}V_{l}^{n}\rangle_{v}\cdot\nabla_{x}K_{l}^{n} +\sum_{p,q=1}^{r}\langle V_{j}^{n}\mathcal{Q}(V_{p}^{n},V_{q}^{n})\rangle_{v}K_{p}^{n}K_{q}^{n},
\end{align}
where the simplification of the collision term relies crucially on its bilinearity and locality in the physical space. Perform a QR decomposition of $K_j^{(1)}= \sum_{i=1}^{r} X_i^{n+1} S_{ij}^{(1)}$ to obtain $X_i^{n+1}$ and $S_{ij}^{(1)}$.

\item \textbf{S-step}: Solve for $S_{ij}^{(2)}$ using the equation
\begin{align} \label{S-scheme}
     \frac{S_{ij}^{(2)}-S_{ij}^{(1)}}{\Delta t}=\sum_{k,l=1}^{r}\langle vV_{j}^{n}V_{l}^{n}\rangle_{v}\cdot\langle X_{i}^{n+1}\nabla_{x}X_{k}^{n+1}\rangle_{x}S_{kl}^{(1)}-\!\!\!\sum_{k,l,p,q=1}^{r}\langle X_{i}^{n+1}X_{k}^{n+1}X_{l}^{n+1}\rangle_{x}\langle V_{j}^{n}\mathcal{Q}(V_{p}^{n},V_{q}^{n})\rangle_{v}S_{kp}^{(1)}S_{lq}^{(1)}.
\end{align}
  \item \textbf{L-step:}  
Define \( L_i^{(2)} = \sum_{j=1}^r S_{ij}^{(2)} V_j^n \), and solve for $L_i^{n+1}$ using the equation
    \begin{align} \label{L-scheme}
     \frac{L_{i}^{n+1}-L_{i}^{(2)}}{\Delta t}=-\sum_{k=1}^{r}v\cdot\langle X_{i}^{n+1}\nabla_{x}X_{k}^{n+1}\rangle_{x}L_{k}^{(2)}+\sum_{k,l,p,q=1}^{r}\langle X_{i}^{n+1}X_{k}^{n+1}X_{l}^{n+1}\rangle_{x}\mathcal{Q}(V_{p}^{n},V_{q}^{n})S_{kp}^{(2)}S_{lq}^{(2)}.
\end{align}
Perform a QR decomposition of $L_i^{n+1}= \sum_{j=1}^{r} S_{ij}^{n+1}V_{j}^{n+1}$ to obtain $S_{ij}^{n+1}$ and $V_{j}^{n+1}$.
\end{itemize}

Note that in the schemes \eqref{K-scheme}, \eqref{S-scheme}, \eqref{L-scheme} above, terms related to the collision operator all involve the form $\mathcal{Q}(V_p^n,V_q^n)$, $p,q=1,\dots,r$. Hence, the collision evaluation only needs to be performed once per time step. This significantly reduces the computational cost, as the collision term always represents the most expensive part of the simulation when solving the Boltzmann equation. Furthermore, the well-developed fast Fourier spectral methods \cite{MP06, GHHH17} can be readily applied to evaluate $\mathcal{Q}(V_p^n,V_q^n)$. More details can be found in \cite{hu2022adaptive}.


\section{Asymptotic-preserving dynamical low-rank method for the stiff Boltzmann equation}
\label{sec:stiff-Boltz}

Although the above dynamical low-rank method significantly reduces the computational cost for solving the Boltzmann equation, the stiffness of the equation, particularly in the regime of small \( \varepsilon \), remains a challenge. In the full tensor framework, for the Boltzmann equation with small \( \varepsilon \), asymptotic-preserving (AP) schemes (see e.g., the review paper \cite{Jin_Rev}) are specifically designed to address the stiffness in the collision term, allowing the collision operator to be handled explicitly while the time step $\Delta t$ can still be independent of $\varepsilon$. 

One such scheme is proposed in \cite{FJ10}. The key idea is to use the BGK operator $\mathcal{Q}_{\text{BGK}}(f) = \mM - f$ to penalize $\mathcal{Q}(f,f)$, where $\mathcal{M}=\mathcal{M}[f]$ is defined in \eqref{Max}. A first-order IMEX scheme for equation \eqref{boltz} reads as
\begin{equation}
    \frac{f^{n+1} - f^n}{\Delta t} + v \cdot \nabla_x f^n = \frac{\mathcal{Q}(f^n, f^n) - \lambda (\mM^n - f^n)}{\varepsilon} + \frac{\lambda (\mM^{n+1} - f^{n+1})}{\varepsilon},\label{eq:Maxwellian distribution-penalty}
\end{equation}
where \( \lambda \) is some positive constant chosen such that $\mathcal{Q}(f,f)\approx \lambda(\mathcal{M}-f)$ and we have used the shorthand notation $\mathcal{M}^n = \mathcal{M}[f^n]$. It has been shown in \cite{FJ10} that the scheme is AP, i.e., \eqref{eq:Maxwellian distribution-penalty} becomes a first-order discretization to the limiting Euler system \eqref{Euler} as $\varepsilon\rightarrow 0$, provided the initial condition is in equilibrium (i.e., is Maxwellian).

It is tempting to extend the full tensor scheme \eqref{eq:Maxwellian distribution-penalty} to the low-rank framework. However, one immediately encounters the difficulty of handling the Maxwellian $\mathcal{M}$. As defined in \eqref{Max}, $\mathcal{M}$ is not a low-rank separated function in $x$ and $v$. Therefore, when deriving the low-rank scheme, one must deal with the projection of the form
\begin{equation}
\langle V_j,\mathcal{M}(x,v,t)\rangle_v, \quad \langle X_i,\mathcal{M}(x,v,t)\rangle_x,
\end{equation}
whose evaluation would require the same computational complexity as a full tensor method, thus negating the advantage of the low-rank approach.

To avoid handling the Maxwellian while still achieving the AP property, \cite{hu2019stochastic} 
 proposed a first-order  scheme by first splitting the advection operator and the collision operator in the Boltzmann equation, and adding a linear penalty in the collision step:
\begin{subequations}\label{fullap}
    \begin{align}
    &\frac{f^* - f^n}{\Delta t} + v \cdot \nabla_x f^n = 0, \label{fullap_a}\\
    &\frac{f^{n+1} - f^*}{\Delta t} = \frac{\mathcal{Q}(f^*, f^*) + \lambda f^* - \lambda f^{n+1}}{\varepsilon}.\label{fullap_b}
    \end{align}
\end{subequations}

By taking the moments $\int_{\mathbb{R}^{d_v}} \cdot \,\phi(v) \rd{v}$ on both sides of \eqref{fullap_b} and using the notations introduced in Section~\ref{sec:intro}, one obtains
\begin{equation}
\frac{\mathcal{U}^{n+1} - \mathcal{U}^*}{\Delta t} = \frac{ \lambda}{\varepsilon} (\mathcal{U}^* -  \mathcal{U}^{n+1}).
\end{equation}
This implies $\mathcal{U}^{n+1}=\mathcal{U}^*$, hence $\mathcal{M}^{n+1}=\mathcal{M}^*$,
i.e., the Maxwellian remains unchanged during the collision step. Taking the moments of \eqref{fullap_a} results in
\begin{equation} \label{moments1}
\frac{\mathcal{U}^*-\mathcal{U}^n}{\Delta t}+\nabla_x\cdot \int  f^n v\phi(v)\rd{v} =0 \ \Longrightarrow \ \frac{\mathcal{U}^{n+1}-\mathcal{U}^n}{\Delta t}+\nabla_x\cdot \int  f^n v\phi(v)\rd{v} =0.
\end{equation}

To formally see the AP property, we rewrite the scheme \eqref{fullap} equivalently as follows 
\begin{align}
f^{n+1}-\mathcal{M}^{n+1}=&\frac{\varepsilon+\mathcal{D}^*\Delta t}{\varepsilon+\lambda \Delta t} (f^n-\mathcal{M}^n)-\frac{\varepsilon \Delta t}{\varepsilon +\lambda\Delta t}\left[v\cdot\nabla_xf^n+\frac{\mathcal{M}^{n+1}-\mathcal{M}^n}{\Delta t}\right] \nonumber\\
&+\frac{\mathcal{D}^*\Delta t^2}{\varepsilon+\lambda \Delta t}\left[\frac{f^*-f^n}{\Delta t}-\frac{\mathcal{M}^*-\mathcal{M}^n}{\Delta t}\right],
\end{align}
where 
\begin{equation}
\mathcal{D}^*:=\frac{\mathcal{Q}(f^*,f^*)}{f^*-\mathcal{M}^*}+\lambda.
\end{equation}
Now, if we assume all the quantities in the above two brackets are $\mathcal{O}(1)$, then as $\varepsilon \rightarrow 0$, we have
\begin{equation}
f^{n+1}-\mathcal{M}^{n+1}=\frac{\mathcal{D}^*}{\lambda}(f^n-\mathcal{M}^n)+\mathcal{O}(\Delta t).
\end{equation}
Suppose we can choose $\lambda$ such that 
\begin{equation} \label{lambda-assumption}
 \frac{\sup|\mathcal{D}^*|}{\lambda}\leq \alpha <1,
\end{equation}
then 
\begin{equation}
|f^{n+1}-\mathcal{M}^{n+1}|\leq \alpha |f^n-\mathcal{M}^n|+\mathcal{O}(\Delta t).
\end{equation}
Iteratively, this gives
\begin{equation}
|f^{n}-\mathcal{M}^{n}|\leq \alpha^n |f^0-\mathcal{M}^0|+\mathcal{O}(\Delta t).
\end{equation}
Therefore, for arbitrary initial data $f^0$, there holds $f^n-\mathcal{M}^n=\mathcal{O}(\Delta t)$ after some time steps. If the initial data is already at the equilibrium $f^0=\mathcal{M}^0$, then $f^n-\mathcal{M}^n=\mathcal{O}(\Delta t)$ for any $n\geq 1$. In either case, there exists an integer $N>0$ such that \eqref{moments1} becomes
\begin{equation} \label{moments-eqn}
\frac{\mathcal{U}^{n+1}-\mathcal{U}^n}{\Delta t}+\nabla_x\cdot \int  \mathcal{M}^n v\phi(v)\rd{v}+\mathcal{O}(\Delta t) =0, \quad \text{for any }n\geq N,
\end{equation}
which is a first-order discretization to the limiting Euler system \eqref{Euler}. Thus the scheme \eqref{fullap} is AP. 

It's not easy to show the existence of $\lambda$ that satisfies the condition \eqref{lambda-assumption} for the Boltzmann collision operator. In practice, it was found in \cite{FJ10, hu2019stochastic} that $\lambda$ can be roughly estimated as
\begin{equation}
\lambda=\sup_{x,v,t}|\mathcal{Q}^-(f)|, \quad \mathcal{Q}^-(f):=\int_{\mathbb{R}^{d_v}}\int_{S^{d_v-1}} B(|v-v_*|,\cos \theta)f(v_*)\rd{\sigma}\rd{v_*}.
\end{equation}
If $\mathcal{Q}$ in \eqref{fullap_b} is the BGK operator $\mathcal{Q}_{\text{BGK}}(f)$, then it is easy to see that $\lambda>1/2$ would satisfy the condition \eqref{lambda-assumption}. In particular, if we choose $\lambda=1$, a stronger AP property can be obtained. As $\varepsilon\rightarrow 0$, \eqref{fullap_b} would imply $f^{n+1}=\mathcal{M}^*=\mathcal{M}^{n+1}$. This means, regardless of the initial condition, we have $f^n=\mathcal{M}^n$ for $n\geq 1$. Therefore, \eqref{moments1} becomes
\begin{equation}
\frac{\mathcal{U}^{n+1}-\mathcal{U}^n}{\Delta t}+\nabla_x\cdot \int  \mathcal{M}^n v\phi(v)\rd{v}=0, \quad \text{for any }n\geq 1,
\end{equation}
which is again a first-order discretization to the limiting system \eqref{Euler} (note that there is no $\mathcal{O}(\Delta t)$ term as in \eqref{moments-eqn}).

Inspired by this full tensor AP scheme, we derive a low-rank scheme as follows. We first split the Boltzmann equation \eqref{boltz} as
\begin{subequations}\label{boltz-split-con}
    \begin{align}
    &\partial_t f = - v \cdot \nabla_x f := \text{RHS}_1,  \label{boltz-split-con_a}\\
    &\partial_t f = \frac{\mathcal{Q}(f, f) + \lambda f - \lambda f}{\varepsilon} := \text{RHS}_2.
    \label{boltz-split-con_b}
    \end{align}
\end{subequations}
 For \eqref{boltz-split-con_a}, which is non-stiff, we can perform a standard first-order explicit {\bf KSL integrator} (just as in Section~\ref{sec:Dynamical low rank methods for the Boltzmann equation} but without the collision term):
\begin{itemize}
\item \textbf{Step 1: K-step}: 
\begin{equation}\label{Step 1}
\begin{split} & K_{j}^{n}=\sum_{i=1}^{r}X_{i}^{n}S_{ij}^{n},\\
 & \frac{K_{j}^{(1)}-K_{j}^{n}}{\Delta t}=-\sum_{l=1}^{r}\langle vV_{j}^{n}V_{l}^{n}\rangle_{v}\cdot\nabla_{x}K_{l}^{n},\\
 & K_{j}^{(1)}=\sum_{i=1}^{r}X_{i}^{*}S_{ij}^{(1)};
\end{split}
\end{equation}
\item \textbf{Step 2: S-step}: 
\begin{equation}\label{Step 2}
\frac{S_{ij}^{(2)}-S_{ij}^{(1)}}{\Delta t}=\sum_{k,l=1}^{r}\langle vV_{j}^{n}V_{l}^{n}\rangle_{v}\cdot\langle X_{i}^{*}\nabla_{x}X_{k}^{*}\rangle_{x}S_{kl}^{(1)};
\end{equation}
\item \textbf{Step 3: L-step}:
\begin{equation}\label{Step 3}
\begin{split} & L_{i}^{(2)}=\sum_{j=1}^{r}S_{ij}^{(2)}V_{j}^{n},\\
 & \frac{L_{i}^{*}-L_{i}^{(2)}}{\Delta t}=-\sum_{k=1}^{r}v\cdot\langle X_{i}^{*}\nabla_{x}X_{k}^{*}\rangle_{x}L_{k}^{(2)},\\
 & L_{i}^{*}=\sum_{j=1}^{r}S_{ij}^{*}V_{j}^{*}.
\end{split}
\end{equation}
\end{itemize}

The stiffness arises in \eqref{boltz-split-con_b}. Directly applying the {\bf KSL integrator} can introduce instability, as it requires solving a backward-in-time S-step. We observe that the collision step fits the framework described in Section~\ref{sec:general XL integrator}. Our proposed {\bf XL integrator} avoids the backward-in-time S-step and reduces the number of differential equations that need to be solved. A first-order IMEX {\bf XL integrator} for \eqref{boltz-split-con_b} is given as follows: 
\begin{itemize}
 \item \textbf{Step 4: X-step} (for \( j = 1,\cdots, r \)):
\begin{equation}\label{K-IMEX}
\begin{split} & K_{j}^{*}=\sum_{i=1}^{r}X_{i}^{*}S_{ij}^{*},\\
 & \frac{K_{j}^{(4)}-K_{j}^{*}}{\Delta t} = \frac{1}{\varepsilon} \sum_{p,q=1}^{r} \langle V_{j}^{*} \mathcal{Q}(V_{p}^{*}, V_{q}^{*}) \rangle_v K_{p}^{*} K_{q}^{*} 
        + \frac{\lambda}{\varepsilon} K_{j}^{*} - \frac{\lambda}{\varepsilon} K_{j}^{(4)},\\
 & [\hat{X}_1^{n+1}, \dots, \hat{X}_{r+\tilde{r}}^{n+1}] = \text{orth}([X_1^*, \dots, X_r^*, K_{1}^{(4)}, \dots, K_{r}^{(4)}]), \quad \tilde{r}=r;
\end{split}
\end{equation}
\item \textbf{Step 5: L-step} (for \( i = 1,\cdots, r+\tilde{r} \)):
\begin{equation}\label{Step 6XL-IMEX}
    \begin{split} 
        & \hat{L}_{i}^{(4)} = 
        \begin{cases} 
        \sum_{j=1}^{r} {S}_{ij}^{*} V_{j}^{*}, & \text{for } 1 \leq i \leq r, \\
        0, & \text{otherwise}.
        \end{cases}\\
        & \frac{\hat{L}_{i}^{n+1} - \hat{L}_{i}^{(4)}}{\Delta t} = \frac{1}{\varepsilon} \sum_{p,q,k,l=1}^{r} \langle \hat{X}_{i}^{n+1} \hat{X}_{k}^{n+1} \hat{X}_{l}^{n+1} \rangle_x \mathcal{Q}(V_{p}^{*}, V_{q}^{*}){S}_{kp}^{*} {S}_{lq}^{*}
        + \frac{\lambda}{\varepsilon} \hat{L}_{i}^{(4)} - \frac{\lambda}{\varepsilon} \hat{L}_{i}^{n+1}, \\
        & \hat{L}_{i}^{n+1} = \sum_{j=1}^{r+\tilde{r}} \hat{S}_{ij}^{n+1} \hat{V}_{j}^{n+1};
    \end{split}
    \end{equation}
     \item \textbf{Step 6: Truncation step}: 
    To reduce the rank back to \( r \), 
    \begin{itemize}
        \item Compute the SVD of \( \hat{S}^{n+1} = (\hat{S}_{ij}^{n+1})_{1\leq i,j\leq r+\tilde{r}} \) to obtain \( \hat{S}^{n+1} = X_s \Sigma_s V_s^T \), where \( X_s, \Sigma_s, V_s \in \mathbb{R}^{(r+\tilde{r}) \times (r+\tilde{r})} \).
        \item Keep the largest \( r \) singular values in \( \Sigma_s \), and set:
        \[
        S^{n+1} = \Sigma_s(1:r, 1:r), \quad X^{n+1} = \hat{X}^{n+1} X_s(:, 1:r), \quad V^{n+1} = \hat{V}^{n+1} V_s(:, 1:r).
        \]
    \end{itemize}
\end{itemize}



By rewriting the differential equation in \eqref{K-IMEX}, we have
\begin{equation}\label{XL-IMEX-1}
    \sum_{i=1}^{r} X_i^{n+1} S_{ij}^{(4)} 
    = \sum_{i=1}^{r} X_i^{*} S_{ij}^{*} 
    + \frac{\Delta t}{\varepsilon + \lambda \Delta t} \sum_{p,q,k,l=1}^{r} \langle V_{j}^{*} \mathcal{Q}(V_{p}^{*}, V_{q}^{*}) \rangle_v S_{kp}^{*} S_{lq}^{*} X_{k}^{*} X_{l}^{*}.
\end{equation}
Clearly, \( \{X_i^{n+1}\}_{i=1}^{r} \) is a linear combination of \( \{X_i^{*}\}_{i=1}^{r} \) and \( \{X_{k}^{*} X_{l}^{*}\}_{k,l=1}^{r} \). Based on this observation, we propose three variants of Step 4, resulting in an {\bf sXL integrator}:
\begin{itemize}
 \item \textbf{Step 4: Augmentation step}:
    \begin{itemize}
        \item Compute $d_{jkl} = \frac{\Delta t}{\varepsilon + \lambda \Delta t}\sum_{p,q=1}^{r}\langle V_{j}^{*}\mathcal{Q}(V_{p}^{*},V_{q}^{*})\rangle_{v}S_{kp}^{*}S_{lq}^{*}$;
        \item Augment the basis functions \( \{X_i^*\}_{i=1}^r \) by including products \( \{X_{k}^{*} X_{l}^{*}\}_{k,l=1}^{r} \):
        \begin{equation}
        [\hat{X}_1^{n+1}, \dots, \hat{X}_{r+\tilde{r}}^{n+1}] = \text{orth}([X_1^*, \dots, X_r^*, X_1^{*} X_1^{*}, X_1^{*} X_2^{*},\dots]).
    \end{equation}
    Three approaches are proposed for this:
    \begin{itemize}
        \item \textbf{Approach 1:} Augment \( \tilde{r}=r \) additional basis functions selected from \( \{X_{k}^{*} X_{l}^{*}\}_{k,l=1}^{r} \) that have the largest coefficients \( |d_{jkl}| \);
        \item \textbf{Approach 2:} Set a tolerance \(\text{tol}\) and include all pairs \((k,l)\) for which \(\vert d_{jkl} \vert > \text{tol}\); 
        \item \textbf{Approach 3:} 
        Augment all \( \tilde{r}=r^2 \) basis functions  \( \{X_{k}^{*} X_{l}^{*}\}_{k,l=1}^{r} \).  
    \end{itemize}
\end{itemize}
\end{itemize}
Note that approach 3 above is a special case of approach 2 when the tolerance is set to zero. This exactly corresponds to the general {\bf sXL integrator} introduced in Section~\ref{sec:general XL integrator}. While accurate, this approach often increases the rank unnecessarily, as demonstrated in numerical examples.

\begin{remark} The XL integrator can be interpreted in the following way: First, perform a K-step to predict the $x$-dependent basis functions. Then perform an L-step to update the coefficients in that basis, which are $v$-dependent functions. This is in contrast to the BUG \cite{ceruti2022unconventional} integrator, where both the $x$- and $v$-dependent basis are predicted and an additional S-step is performed to update the coefficients (which are $r \times r$ time-dependent matrices). Nevertheless, we borrow from \cite{ceruti2022rank} the augmentation procedure that is used to obtain $\hat{X}$. This is important for the XL integrator as well, as otherwise the initial value for the L-step is not exactly represented and this can result in accuracy and stability issues (see, e.g., the simulation with the BUG integrator in \cite{Einkemmer2023}). We note that due to the smaller number of steps the XL integrator has lower computational cost than the BUG integrator, though the asymptotic complexity is the same for both methods.
\end{remark}

\begin{remark}
The separable property that allows us to obtain the sXL integrator is not only satisfied by the Boltzmann collision operator. Indeed, many kinetic operators have a similar structure. We mention some of them here and defer their numerical exploration to future work. Assume that $f(x,v,t)=\sum_{j=1}^r K_j(x,t)V_j(v,t)$, then
\begin{itemize}
\item the advection operator:
\[
v\partial_{x}f=\sum_{j=1}^r\partial_{x}K_{j}vV_{j};
\]
\item the Fokker--Planck type operator:
\[
\partial_{v}(a(v)f)+\partial_{vv}(b(v)f)=\sum_{j=1}^rK_{j}\left[\partial_{v}(a(v)V_j)+\partial_{vv}(b(v)V_j)\right];
\]
\item the special BGK operator, where the Maxwellian term $M$ is low-rank (such as the model considered in Section~\ref{sec:BGK operator and proof of AP property}):
\[
M-f=\sum_{j=1}^rK_{j}^{M}V_{j}^{M}-\sum_{j=1}^rK_{j}V_{j}.
\]
\end{itemize}
\end{remark}

\section{The AP property for a special BGK equation}
\label{sec:BGK operator and proof of AP property}

To mathematically prove that the AP property holds for the proposed {\bf XL} and {\bf sXL} integrators for the nonlinear Boltzmann equation is challenging. Even replacing the Boltzmann operator by the BGK operator is still a difficult problem. The reason for this is that even though the BGK operator is significantly simpler it has the added complication (not found in the Boltzmann collision operator) that no general low-rank representation of the Maxwellian is available (see, e.g., \cite{chen1998lattice, einkemmer2019low, Dektor2024}). 

We thus focus on a special 1D1V ($d_x=d_v=1$) BGK equation:
\begin{equation} \label{BGK1}
\partial_t f + v \partial_x f = \frac{1}{\varepsilon}(M[f] - f),
\end{equation}
where the (simplified) Maxwellian $M[f]$ is given by
\begin{equation}
M[f] := \rho \frac{e^{-v^2 / 2}}{\sqrt{2\pi}} \left(1 + vu + \frac{1}{2}(v^2 - 1)u^2 \right), \quad \rho=\int_{\mathbb{R}} f\rd{v}, \quad u=\frac{1}{\rho}  \int_{\mathbb{R}} fv\rd{v}.
\end{equation}
It is easy to check that
\begin{equation}
\int_{\mathbb{R}}f \varphi(v)\rd{v}=\int_{\mathbb{R}}M\varphi(v)\rd{v}=(\rho,\rho u)^T:=U, \quad \varphi(v):=(1,v)^T.
\end{equation}
Now, taking the moments $\int_{\mathbb{R}}\cdot \ \varphi(v)\rd{v}$ on both sides of \eqref{BGK1} yields
\begin{equation} \label{local-con1}
\partial_t U+\partial_x \int_{\mathbb{R}} fv \varphi(v)\rd{v}=0.
\end{equation}
As $\varepsilon \rightarrow 0$, \eqref{BGK1} implies $f\rightarrow M[f]$, then \eqref{local-con1} becomes a closed hyperbolic system:
\begin{align} \label{hyper}
\partial_t U+\partial_x \int_{\mathbb{R}} M[f]v \varphi(v)\rd{v}=0  \ \Longleftrightarrow
\left\{
\begin{array}{l}
\displaystyle \partial_t \rho+\partial_x (\rho u)=0,\\[8pt]
\displaystyle \partial_t (\rho u)+\partial_x (\rho u+\rho u^2)=0.
\end{array}
\right.
\end{align}

Note that the Maxwellian function $M[f]$ can be obtained from the general Maxwellian $\mathcal{M}[f]$ \eqref{Max} by assuming constant temperature and small bulk velocity (it is a good approximation in the weakly compressible regime \cite{einkemmer2019low}). $M[f]$ is a clear low-rank (rank-3) function, making the following analysis of the low-rank algorithm tractable. 

The {\bf XL} and {\bf sXL} integrators introduced in Section~\ref{sec:stiff-Boltz}, adapted to the special BGK equation \eqref{BGK1}, is outlined as follows. In particular, this means $\mathcal{Q} = M[f] - f$ and $\lambda = 1$ in \eqref{boltz-split-con_b}. To avoid redundancy, we only present Steps 4-5. Steps 1-3 and Step 6 remain the same as those in Section~\ref{sec:stiff-Boltz}.
\begin{itemize}
 \item \textbf{Step 4 } (for \( j = 1,\cdots, r \)):
\begin{itemize}
    \item Define
    \begin{equation}
K_{j}^{*}=\sum_{i=1}^{r}X_{i}^{*}S_{ij}^{*},
    \end{equation}
    and compute $M^*$ as
    \begin{equation} \label{computeM}
     \rho^* = \sum_{j=1}^{r} K_j^*\langle V_j^*\rangle_v, \quad 
    u^* = \frac{1}{\rho^*} \sum_{j=1}^{r} K_j^*\langle vV_j^*\rangle_v, \quad
        M^{*} = \rho^{*} \frac{e^{-v^{2}/2}}{\sqrt{2\pi}} \left( 1 + v u^{*} + \frac{1}{2} \left( v^{2} - 1 \right) (u^{*})^{2} \right);
    \end{equation}
    \item \textbf{XL integrator}
    \begin{equation}
    \begin{split} 
     & \frac{K_{j}^{(4)}-K_{j}^{*}}{\Delta t} = \frac{1}{\varepsilon} \langle M^{*}V_{j}^{*}  \rangle_v 
            - \frac{1}{\varepsilon} K_{j}^{(4)},\\
     & [\hat{X}_1^{n+1}, \dots, \hat{X}_{2r}^{n+1}] = \text{orth}([X_1^*, \dots, X_r^*, K_{1}^{(4)}, \dots, K_{r}^{(4)}]);
    \end{split}
    \end{equation}
    \item \textbf{sXL integrator}
    \begin{equation}
    [\hat{X}_1^{n+1}, \dots, \hat{X}_{2r}^{n+1}] = \text{orth}([X_1^*, \dots, X_r^*, \langle M^{*}V_{1}^{*}  \rangle_v , \dots, \langle M^{*}V_{r}^{*}  \rangle_v ]).
    \end{equation}
\end{itemize}
\item \textbf{Step 5 } (for \( i = 1,\cdots, 2r \)):
\begin{equation}
\begin{split}
 & \hat{L}_{i}^{(4)} = 
        \begin{cases} 
        \sum_{j=1}^{r} {S}_{ij}^{*} V_{j}^{*}, & \text{for } 1 \leq i \leq r, \\
        0, & \text{otherwise}.
        \end{cases}\\
 & \frac{\hat{L}_{i}^{n+1} - \hat{L}_{i}^{(4)}}{\Delta t} = \frac{1}{\varepsilon} \langle \hat{X}_{i}^{n+1}M^{*} \rangle_{x} - \frac{1}{\varepsilon} \hat{L}_{i}^{n+1},\\
  &\hat{L}_{i}^{n+1}=\sum_{j=1}^{2r}\hat{S}_{ij}^{n+1}\hat{V}_{j}^{n+1}.
\end{split}
\end{equation}
\end{itemize}

We have the following theorem for the AP property:
\begin{theorem} [\textbf{AP property}]
    Assume that the initial data is at the equilibrium $f^0=\sum_{j=1}^{r} K_j^0 V_j^0=M^0$, and rank $r\geq 6$ is used in the above low-rank algorithm. Then as $\varepsilon \rightarrow 0$, we have for any integer $n\geq 0$,
    \begin{equation}
    f^{n}=\sum_{j=1}^{r} K_j^{n} V_j^{n}=M^{n}, 
    \end{equation}
    and 
\begin{equation} \label{UU}
    \frac{U^{n+1}-U^n}{\Delta t} = -\partial_x \left( \begin{matrix} \rho^n u^n \\ \rho^n+\rho^n (u^n)^2 \end{matrix} \right) +  \mathcal{O}(\Delta t),
\end{equation}
where $U^n=\int_{\mathbb{R}}f^n\varphi(v)\rd{v}$. That is, the low-rank algorithm becomes a first-order discretization to the limiting system (\ref{hyper}), hence is AP.
\end{theorem}
\begin{proof}
It suffices to assume $f^n=M^n$, and prove $f^{n+1}=M^{n+1}$ and \eqref{UU}.

First, observe that 
\begin{align*}
    f^n = M^n = \rho^n \frac{e^{-v^{2}/2}}{\sqrt{2\pi}} \left(1 + vu^n + \frac{1}{2}(v^{2}-1)(u^n)^{2}\right)=\sum_{j=1}^{r} K_j^n V_j^n.
\end{align*}
By applying a QR decomposition, we can express the basis functions $\{V_j^n\}_{j=1}^r$ as
\begin{align*}
    V_1^{n} = \alpha_1 e^{-v^2/2}, \
    V_2^{n} = 2\alpha_2 v e^{-v^2/2}, \
    V_3^{n} = \alpha_3 \left(4v^2 - 2\right) e^{-v^2/2}, \
    V_4^{n} = \alpha_4\left(8v^3 - 12v\right) e^{-v^2/2}, \dots,
\end{align*}
and \( K_j^{n} = 0 \) for $j\geq 4$. The coefficients \( \alpha_1, \alpha_2, \ldots, \alpha_r \) are determined to ensure that the functions \( V_1^n, V_2^n, \ldots, V_r^n \) remain orthonormal. Utilizing the recurrence relations associated with Hermite polynomials, we derive the following recurrence relations:
\begin{equation*}
\left\{\begin{split}
& vV_1^n = C_2V_2^n, \\
& vV_j^n ={C_{j+1}}V_{j+1}^n + {D_{j-1}}V_{j-1}^n  \ 
\text{ for } j \geq 2,
\end{split}\right.\end{equation*}
where $C_{j+1} = \frac{\alpha_j}{2\alpha_{j+1}}\text{ for } j \geq 1$ and $D_{j-1} = \frac{\alpha_j(j-1)}{\alpha_{j-1}} \text{ for } j \geq 2$. By defining $D_0 = 0$ and $V_0^n = 0$, we have a unified recurrence relation as
\begin{equation}\label{recurrence}
    vV_j^n ={C_{j+1}}V_{j+1}^n + {D_{j-1}}V_{j-1}^n \ \text{ for } j \geq 1.
\end{equation}

Now, let's dive into each step.

Step 1:
\begin{equation}\label{apstep1-1}
    \frac{K_{j}^{(1)}-K_{j}^{n}}{\Delta t}=-\sum_{l=1}^{r}\langle vV_{j}^{n}V_{l}^{n}\rangle_{v}\partial_x K_{l}^{n} = -\sum_{l=1}^{3}\langle vV_{j}^{n}V_{l}^{n}\rangle_{v}\partial_x K_{l}^{n}.
\end{equation}
We also have
\begin{align*}
    &\rho^n u^n = \int f^n v \rd{v} = \sum_{i=1}^{3} K_{i}^n \langle v V_{i}^n \rangle_{v}, \quad 
    \rho^n+\rho^n (u^n)^2 = \int  f^nv^2  \rd{v}= \sum_{i=1}^{3} K_{i}^n \langle v^2 V_{i}^n \rangle_{v}.
\end{align*}
After Step 1, the macroscopic quantities can be obtained as follows:
\begin{equation*}
    U^{(1)} := \sum_{j=1}^{r} K_{j}^{(1)}\langle \varphi(v) V_{j}^{n}\rangle_{v}.
\end{equation*}
Taking $\sum_{j=1}^r \cdot \ \langle \varphi(v)V_j^n\rangle_v$ on both sides of \eqref{apstep1-1} yields
\begin{equation}\label{rho1}
    \begin{aligned}
    \frac{U^{(1)} - U^{n}}{\Delta t} =& \sum_{j=1}^{r} \frac{K_{j}^{(1)}- K_{j}^{n}}{\Delta t}\langle \varphi(v)V_{j}^{n}\rangle_{v}\\
    =&- \sum_{j=1}^{r}\sum_{l=1}^{3}\langle vV_{j}^{n}V_{l}^{n}\rangle_{v}\partial_x K_{l}^{n}\langle \varphi(v)V_{j}^{n}\rangle_{v},\\
    =& -\sum_{j=1}^{r}\sum_{l=1}^{3}\langle V_{j}^{n}({C_{l+1}}V_{l+1}^n + {D_{l-1}}V_{l-1}^n)\rangle_{v}\partial_x K_{l}^{n}\langle \varphi(v)V_{j}^{n}\rangle_{v}\\
    =& -\sum_{l=1}^{3}\langle V_{l+1}^{n}{C_{l+1}}V_{l+1}^n \rangle_{v}\partial_x K_{l}^{n}\langle \varphi(v)V_{l+1}^{n}\rangle_{v}-\sum_{l=1}^{3}\langle V_{l-1}^{n} {D_{l-1}}V_{l-1}^n\rangle_{v}\partial_x K_{l}^{n}\langle \varphi(v)V_{l-1}^{n}\rangle_{v}\\
    =& -\sum_{l=1}^{3} \partial_x K_{l}^{n}\langle \varphi(v)(C_{l+1}V_{l+1}^{n}+D_{l-1}V_{l-1}^{n})\rangle_{v}\\
    =& -\sum_{l=1}^{3} \partial_x K_{l}^{n}\langle \varphi(v)vV_l^n\rangle_{v}\\
    =& -\partial_x \left( \begin{matrix} \rho^n u^n \\ \rho^n+\rho^n (u^n)^2 \end{matrix} \right).
    \end{aligned}
\end{equation}

Step 2 and Step 3: 
We can multiply \eqref{Step 2} by \( \sum_{j=1}^r \cdot \ \varphi(v)V_j^n \), multiply \eqref{Step 3} by \( \varphi(v) \), and sum them to obtain
\begin{equation*}
\begin{aligned}
&\sum_{j=1}^{r}\frac{S_{ij}^{(2)}\varphi(v)V_j^n-S_{ij}^{(1)}\varphi(v)V_j^n}{\Delta t} + \frac{\varphi(v)L_{i}^{*}-\varphi(v)L_{i}^{(2)}}{\Delta t}\\ =&\sum_{j, k,l=1}^{r}\langle vV_{j}^{n}V_{l}^{n}\rangle_{v}\langle X_{i}^{*}\partial_{x}X_{k}^{*}\rangle_{x}S_{kl}^{(1)}\varphi(v)V_j^n-\sum_{{k=1}}^{r}v\langle X_{i}^{*}\partial_{x}X_{k}^{*}\rangle_{x}\varphi(v)L_{k}^{(2)}.
\end{aligned}
\end{equation*}
By denoting $L_i^{(1)} = \sum_{j=1}^{r}S_{ij}^{(1)}V_j^n$, and using $L_i^{(2)} = \sum_{j=1}^{r}S_{ij}^{(2)}V_j^n$, we have
\begin{equation*}
\begin{aligned}
    \varphi(v){L_{i}^{*}}&=\varphi(v)L_{i}^{(1)} + {\Delta t}\sum_{j, k,l=1}^{r}\langle X_{i}^{*}\partial_{x}X_{k}^{*}\rangle_{x}\langle vV_{j}^{n}V_{l}^{n}\rangle_{v}S_{kl}^{(1)}\varphi(v)V_j^n-{\Delta t}\sum_{k,l=1}^{r}\langle X_{i}^{*}\partial_{x}X_{k}^{*}\rangle_{x}S_{kl}^{(2)}v\varphi(v){V_l^n}\\
    &=\varphi(v)L_{i}^{(1)} + {\Delta t}\sum_{k=1}^{r}\langle X_{i}^{*}\partial_{x}X_{k}^{*}\rangle_{x}\left(\sum_{j, l=1}^{r}\langle vV_{j}^{n}V_{l}^{n}\rangle_{v}S_{kl}^{(1)}\varphi(v)V_j^n-\sum_{l=1}^{r}S_{kl}^{(2)}v\varphi(v){V_l^n}\right).
    \end{aligned}
    \end{equation*} 
 We then take the integral of the above equation with respect to \( v \) to obtain
\begin{equation}\label{L1BGK}
\begin{aligned}
    \langle{\varphi(v)L_{i}^{*}}\rangle_{v} -  \langle{\varphi(v)L_{i}^{(1)}}\rangle_{v}&= {\Delta t}\sum_{k=1}^{r}\langle X_{i}^{*}\partial_{x}X_{k}^{*}\rangle_{x}\left\langle\sum_{j, l=1}^{r}\langle vV_{j}^{n}V_{l}^{n}\rangle_{v}S_{kl}^{(1)}\varphi(v)V_j^n-\sum_{l=1}^{r}S_{kl}^{(2)}v\varphi(v){V_l^n}\right\rangle_{v}.
    \end{aligned}
    \end{equation} 
We now focus on the big bracket term in \eqref{L1BGK} and prove that 
\begin{equation} \label{bigbracket}
\left\langle\sum_{j, l=1}^{r}\langle vV_{j}^{n}V_{l}^{n}\rangle_{v}S_{kl}^{(1)}\varphi(v)V_j^n-\sum_{l=1}^{r}S_{kl}^{(2)}v\varphi(v){V_l^n}\right\rangle_{v}=\mathcal{O}(\Delta t),
\end{equation}
then
\begin{equation} \label{L-estimate}
\langle{\varphi(v)L_{i}^{*}}\rangle_{v} -  \langle{\varphi(v)L_{i}^{(1)}}\rangle_{v}=\mathcal{O}(\Delta t^2).
\end{equation}
To see this, we first write \eqref{bigbracket} as
\begin{equation} \label{bigbracket1}
\begin{aligned}
&\left\langle\sum_{j, l=1}^{r}\langle vV_{j}^{n}V_{l}^{n}\rangle_{v}S_{kl}^{(1)}\varphi(v)V_j^n-\sum_{l=1}^{r}S_{kl}^{(2)}v\varphi(v){V_l^n}\right\rangle_{v}\\
=&\left\langle \sum_{l=1}^{r-1}\left(\sum_{j=1}^r \langle vV_{j}^{n}V_{l}^{n}\rangle_{v}\varphi(v)V_j^nS_{kl}^{(1)}-v\varphi(v){V_l^n}S_{kl}^{(2)}\right) \right\rangle_{v} +\left\langle\sum_{j=1}^{r}\langle vV_{j}^{n}V_{r}^{n}\rangle_{v} \varphi(v)V_j^n S_{kr}^{(1)}-v\varphi(v)V_{r}^{n}S_{kr}^{(2)}\right\rangle_{v}.
\end{aligned}
\end{equation}
For the first term on the right hand side of \eqref{bigbracket1}, using \eqref{recurrence}, we have for $1\leq l \leq r-1$, 
    \begin{equation}
    \begin{aligned}
        \sum_{j=1}^{r}\langle vV_{j}^{n}V_{l}^{n}\rangle_{v} \varphi(v)V_j^n &=  \sum_{j=1}^{r}\langle V_{j}^{n}\left( {C_{l+1}}V_{l+1}^n + {D_{l-1}}V_{l-1}^n \right)\rangle_{v} \varphi(v)V_j^n\\
        &=  \left(C_{l+1}V_{l+1}^n+ D_{l-1}V_{l-1}^n\right) \varphi(v)\\
        &=  v \varphi(v) V_l^{n}.
    \end{aligned}
    \end{equation}
Thus,
\begin{equation}
\left\langle \sum_{l=1}^{r-1}\left(\sum_{j=1}^r \langle vV_{j}^{n}V_{l}^{n}\rangle_{v}\varphi(v)V_j^nS_{kl}^{(1)}-v\varphi(v){V_l^n}S_{kl}^{(2)}\right) \right\rangle_{v} =\sum_{l=1}^{r-1} \left\langle v\varphi(v)V_l^n\right\rangle_v \left(S_{kl}^{(1)}-S_{kl}^{(2)}\right).
\end{equation}
Using \eqref{Step 2}, we have 
\begin{equation}\label{SijOt}
     S_{ij}^{(1)} - S_{ij}^{(2)} = \mathcal{O}(\Delta t)  \ \text{ for any } 1 \leq i,j \leq r.
\end{equation}
For the second term on the right hand side of \eqref{bigbracket1}, we are going to show 
\begin{equation}\label{S1S2=0}
    S_{ij}^{(1)} =S_{ij}^{(2)}= 0 \ \text{ for } 1\leq i \leq r, ~6 \leq j \leq r,
\end{equation}
which follows $S_{kr}^{(1)} =S_{kr}^{(2)}= 0$ for $k=1, \cdots r$. For $j \geq 5$, utilizing \eqref{apstep1-1} and \eqref{recurrence}, we have
\begin{equation} \label{Kr50}
\frac{K_{j}^{(1)}}{\Delta t} = -\sum_{l=1}^{3} \langle v V_{j}^{n} V_{l}^{n} \rangle_{v} \partial_{x} K_{l}^{n} = 0.
\end{equation}
By taking the inner product of the above equation with \( X_i^{*} \), we have
\begin{equation*}
    S_{ij}^{(1)} = -{\Delta t}\sum_{l=1}^{3} \langle v V_{j}^{n} V_{l}^{n} \rangle_{v} \langle X_k^{*}\partial_{x} K_{l}^{n}  \rangle_{x}=0, \ \text{ for } j \geq 5.
\end{equation*} 
For $j \geq 6$,
\begin{equation*}
\begin{aligned}
    {S_{ij}^{(2)}}&=S_{ij}^{(1)} + {\Delta t}\sum_{k,l=1}^{r}\langle vV_{j}^{n}V_{l}^{n}\rangle_{v}\langle X_{i}^{*}\partial_x X_{k}^{*}\rangle_{x}S_{kl}^{(1)}\\
    &=S_{ij}^{(1)} + {\Delta t}\sum_{k=1}^{r}\sum_{l=1}^{4}\langle vV_{j}^{n}V_{l}^{n}\rangle_{v}\langle X_{i}^{*}\partial_x X_{k}^{*}\rangle_{x}S_{kl}^{(1)}\\
    &=0.
\end{aligned}
\end{equation*}
Combining everything, we have established \eqref{bigbracket}.

Now using \eqref{L-estimate} and defining $U^{*} := \sum_{i=1}^{r} X_{i}^{*}\langle  \varphi(v)L_{i}^{*}\rangle_{v}$,
we have
\begin{align}
    &\frac{U^{*} - U^{(1)}}{\Delta t}=\frac{1}{\Delta t}\sum_{i=1}^{r} X_{i}^{*} \langle {\varphi(v) L_{i}^{*}} - \varphi(v)L_{i}^{(1)} \rangle_{v}={\mathcal{O}(\Delta t).}
\end{align}
Combining the above with \eqref{rho1}, we have
\begin{align}\label{rho2}
    &\frac{U^{*} - U^{n}}{\Delta t}=-\partial_x\left( \begin{matrix} \rho^n u^n \\ \rho^n+\rho^n (u^n)^2 \end{matrix} \right)+{\mathcal{O}(\Delta t).}
\end{align}

At the end of Step 3, we perform a QR decomposition of 
\[
L_{i}^{*} = \sum_{j=1}^{r} S_{ij}^{*} V_{j}^{*}.
\] 
From \eqref{S1S2=0}, we have
\[
L_{i}^{(2)} = \sum_{j=1}^{r} S_{ij}^{(2)} V_{j}^{n} = \sum_{j=1}^{5} S_{ij}^{(2)} V_{j}^{n}.
\]
Thus,
\[
L_{i}^{*} = L_{i}^{(2)} - \Delta t \sum_{k=1}^{r} \langle X_{i}^{*} \, \partial_x X_{k}^{*} \rangle_{x} \, v L_{k}^{(2)} = \sum_{j=1}^{5} S_{ij}^{(2)} V_{j}^{n} - \Delta t \sum_{k=1}^{r} \langle X_{i}^{*} \, \partial_x X_{k}^{*} \rangle_{x} \, v \sum_{j=1}^{5} S_{kj}^{(2)} V_{j}^{n}.
\]
This implies that \( L_i^{*} \in \text{span}\left\{ V_j^{n} \right\}_{j=1,2,\ldots, r} \) for \( r \geq 6 \) and $V_j^{*} = V_j^{n}$ for $j=1,2,\ldots, r$ with proper QR decomposition.

Step 4: For $M^{*}$ defined in \eqref{computeM}, it is clear that \( M^{*} \in \text{span}\left\{ V_j^{*} \right\}_{j=1,2,3} \). For \textbf{XL integrator}, when $\varepsilon \rightarrow 0$, we have
\begin{equation}
    K_j^{(4)} = \langle M^{*}V_j^{*}\rangle_{v},\  \text{ for } 1 \leq j \leq r,
\end{equation}
so
\begin{equation}
     [\hat{X}_1^{n+1}, \dots, \hat{X}_{2r}^{n+1}] = [X_{1}^{*}, \dots, X_{r}^{*}, \langle M^{*}V_1^{*}\rangle_{v}, \dots, \langle M^{*}V_r^{*}\rangle_{v}],
\end{equation}
which is the same as \textbf{sXL integrator}.

After step 4, we have
\begin{equation}
    f^{(4)} := \sum_{i=1}^{2r}\hat{X}_i^{n+1}\hat{L}_i^{(4)} = \sum_{i=1}^{r}{X}_i^{*}\hat{L}_i^{(4)} = \sum_{i,j=1}^{r}{X}_i^{*}S_{ij}^{*}V_j^{*} = f^{*}.
\end{equation}
where the second equality holds because \( [X_{1}^{*}, \dots, X_{r}^{*}] \) is already orthonormal. It follows that the Maxwellian \( M \), which depends on the moments of \( f \), also remains unchanged. So we have
\begin{equation}\label{M n+1 star star}
    M^{*} \in \text{span}\left\{ \hat{X}_i^{n+1}\right\}_{i=1,\cdots, 2r}.
\end{equation}

Step 5:  
When $\varepsilon \rightarrow 0$, we have
$$\hat{L}_i^{n+1} = \langle\hat{X}_i^{n+1}M^{*}\rangle_{x}, \ \text{ for } 1 \leq i \leq 2r.$$
Then
\begin{equation}\label{xl-fnew}
    \begin{aligned}
        \hat{f}^{n+1}=\sum_{i=1}^{{2r}} \hat{X}_i^{n+1}\hat{L}_i^{n+1} =  \sum_{i=1}^{{2r}} \hat{X}_i^{n+1}\langle \hat{X}_i^{n+1}M^{*}\rangle_{x} = M^{*}.
    \end{aligned}
\end{equation}
Hence, 
\begin{align}
    &\hat{U}^{n+1} =\sum_{i=1}^{2r} \hat{X}_i^{n+1}\langle \varphi(v)\hat{L}_i^{n+1}\rangle_{v}= \langle \varphi(v)M^{*}\rangle_{v}=U^{*}.\label{xl-rhonew}
\end{align}

Step 6:
Since \( M^* \) lies within the span of \( \{V_j^{*}\}_{j=1,2,3} \), $\hat{f}^{n+1}$ and \( \hat{L}_i^{n+1} \) also lie within the span of \( \{V_j^{*}\}_{j=1,2,3} \). By performing QR decompositions on \( \hat{L}_i^{n+1} \) to obtain \( \hat{S}^{n+1} \) and \( \hat{V}^{n+1} \) and applying SVD to \( \hat{S}^{n+1} = X_s \Sigma_s V_s^T \) results in:
\[
\Sigma_s = \begin{bmatrix} \Sigma_{3 \times 3} & 0 \\ 0 & 0 \end{bmatrix}.
\]
The rank of \( \hat{S}^{n+1} \) remains 3, and truncation at this step back to rank $r\geq 6$ does not result in any change of $\hat{f}^{n+1}$. Hence,
\begin{equation}\label{equalfM}
    {f}^{n+1} = \hat{f}^{n+1} = M^{*}, \quad {U}^{n+1} = \hat{U}^{n+1} = U^{*}.
\end{equation}

Combining the above with equation \eqref{rho2}, we obtain
\begin{align}
    &\frac{U^{n+1} - U^{n}}{\Delta t}=-\partial_x\left( \begin{matrix} \rho^n u^n \\ \rho^n+\rho^n (u^n)^2 \end{matrix} \right) + \mathcal{O}(\Delta t).\label{rho3Wronglabel}
\end{align}
Furthermore, \eqref{equalfM} implies \( M^{n+1} = M^{*} \), and thus
\[
f^{n+1} = M^{n+1}.
\]
\end{proof}

\section{Numerical results}
\label{sec:experiments}

In this section, we present several numerical experiments to evaluate the performance of the proposed XL integrator (DLR-XL) alongside its variants, the sXL integrator with the various strategies for adding basis functions as outlined in Section \ref{sec:stiff-Boltz}. We denote the sXL integrator with approach 1 as DLR-sXL-1 and with approach 2 as DLR-sXL-2. The full tensor method is implemented following the scheme \eqref{fullap_a}-\eqref{fullap_b}. Unless otherwise specified, the following assumptions and configurations are used:
\begin{itemize}
    \item \textbf{Spatial domain}: The spatial variable \( x \) lies in the interval \([0,1]\) ($d_x=1$). For the first example, we impose periodic boundary conditions, while Neumann boundary conditions are applied in the second example. Spatial discretization is performed using the second-order MUSCL scheme \cite{van1979towards} with a minmod slope limiter. The spatial domain is discretized with 100 grid points.
    
    \item \textbf{Velocity domain}: The velocity variable \( v \) is defined over the domain \([ -L_v, L_v]^2\) ($d_v=2$), where \( L_v = 8.4 \). The Boltzmann collision operator \(\mathcal{Q}\) with constant kernel $B=1/(2\pi)$ is evaluated using the fast spectral method introduced in \cite{MP06}. Each velocity dimension is discretized with 32 points for all examples. The penalty parameter \(\lambda\) is set to 1.1.

    \item \textbf{Rank}: \( r \) denotes the initial rank set for low-rank methods. \( \tilde{r} \) represents the additional rank augmented in Step 4 for low-rank methods, fixed at \( r \) for DLR-XL and DLR-sXL-1, and dependent on the tolerance, \(\text{tol}\), for DLR-sXL-2. \(\text{tol}\) can be chosen based on prior information or tested iteratively from large to small values, progressively increasing the additional rank.     \item \textbf{Macroscopic quantities}: The macroscopic quantities, density $\rho$, bulk velocity $u$, temperature $T$, of full tensor method are computed using \eqref{macro-quantity}.
    For the low-rank methods, they can be computed directly from the low-rank factors:
    \[
    \rho = \sum_{j=1}^{r} K_j\langle V_j\rangle_v, \quad 
    u = \frac{1}{\rho} \sum_{j=1}^{r} K_j\langle vV_j\rangle_v,  \quad 
    T = \frac{1}{d_v\rho} \sum_{j=1}^{r} K_j\langle |v|^2V_j\rangle_v - \frac{1}{d_v}|u|^2.
    \]
    For the bulk velocity \( u \), which is a two-dimensional vector, we will present results only for its first component for brevity.
\end{itemize}

\subsection{Sine-type initial data}

We first evaluate the accuracy of our proposed integrators compared with the full tensor scheme using a continuous sine-type initial condition in both the kinetic and fluid regimes. Specifically, we consider a large \(\varepsilon = 1\) and a very small \(\varepsilon = 10^{-6}\), and analyze the macroscopic quantities. 

The initial data is at the equilibrium:
\begin{equation}
f^0(x, v) = \frac{\rho^0}{2\pi T^0} e^{-\frac{|v-u^0|^2}{2T^0}},
\end{equation}
where
\begin{equation}
\rho^0(x) = \frac{2 + \sin(2\pi x)}{3}, \quad u^0 = (0.2, 0), \quad T^0(x) = \frac{3 + \cos(2\pi x)}{4}.
\end{equation}

\begin{figure}[htbp]
    \centering
     \begin{subfigure}[b]{0.48\textwidth}
        \includegraphics[width=\textwidth]{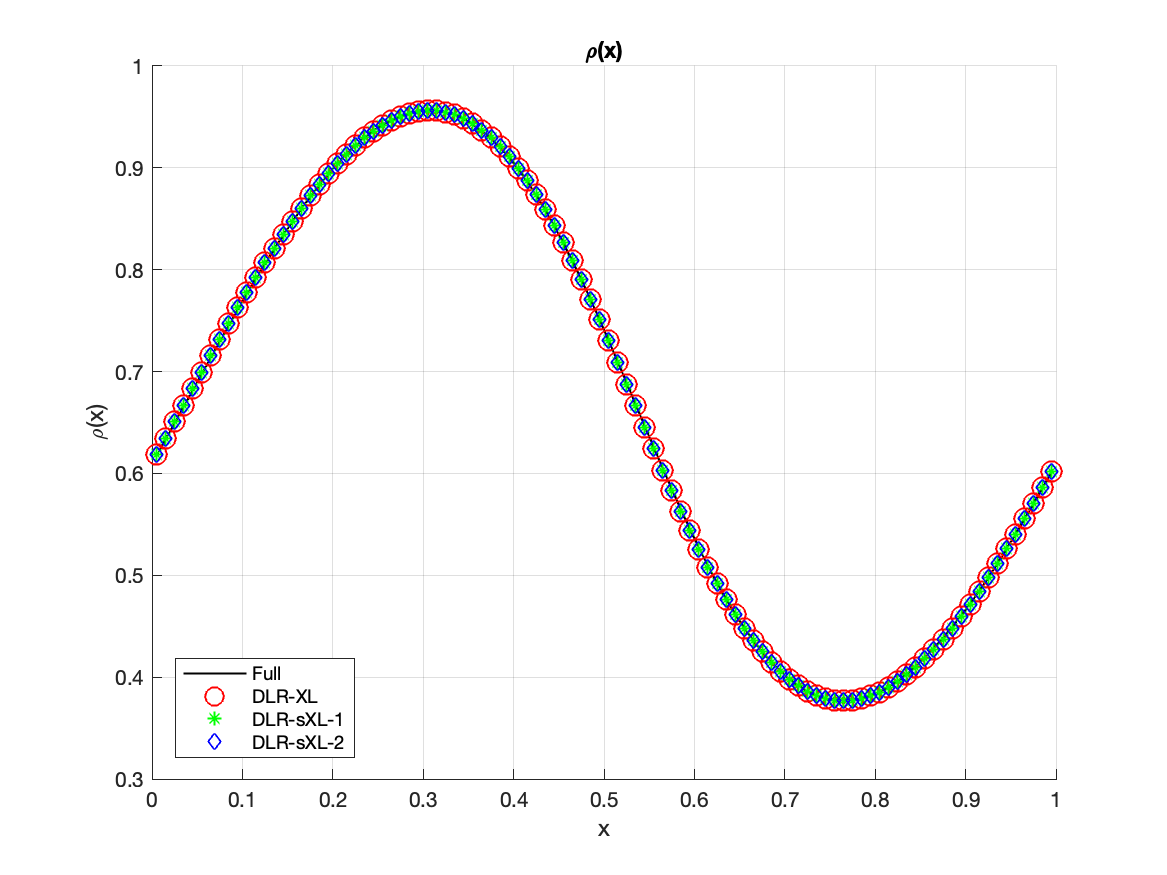}
        \caption{Density \(\rho\).}
    \end{subfigure}
    \begin{subfigure}[b]{0.48\textwidth}
        \includegraphics[width=\textwidth]{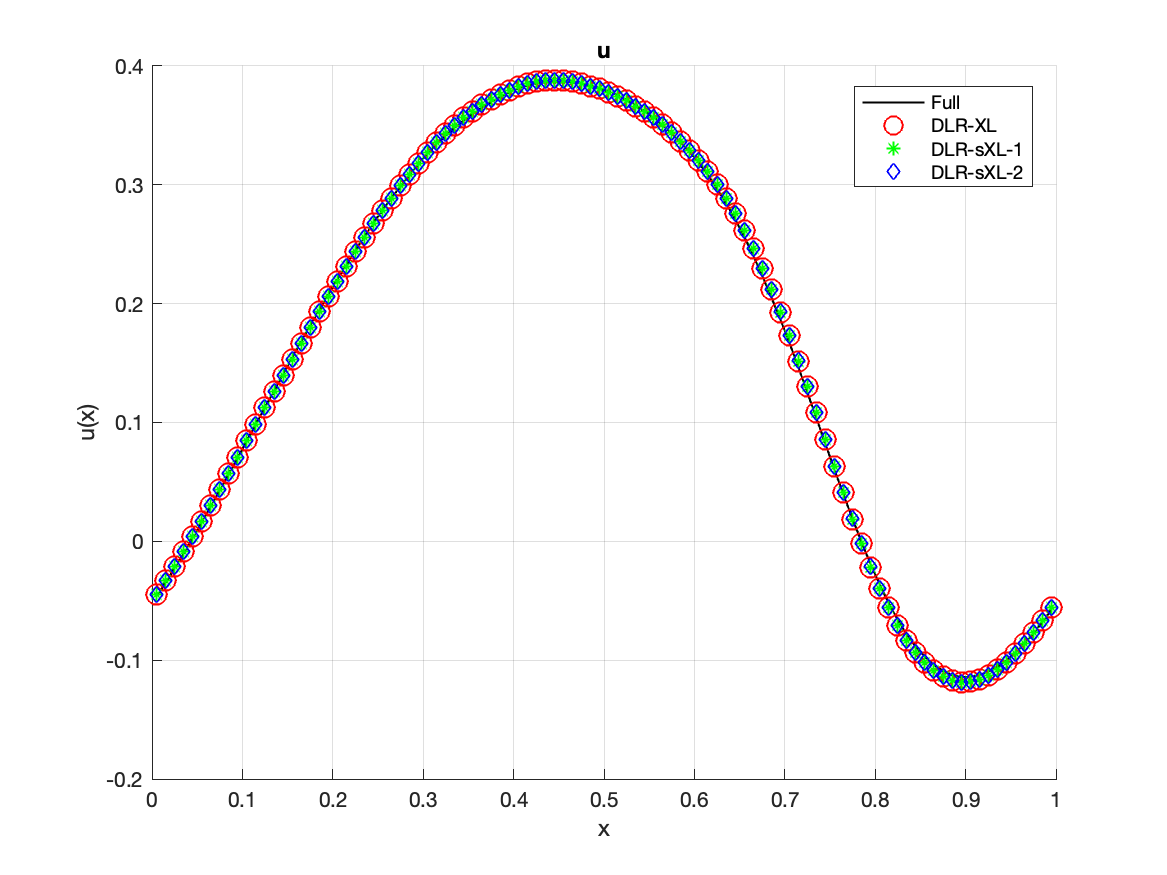}
        \caption{Bulk velocity \(u\).}
    \end{subfigure}
    \begin{subfigure}[b]{0.48\textwidth}
        \includegraphics[width=\textwidth]{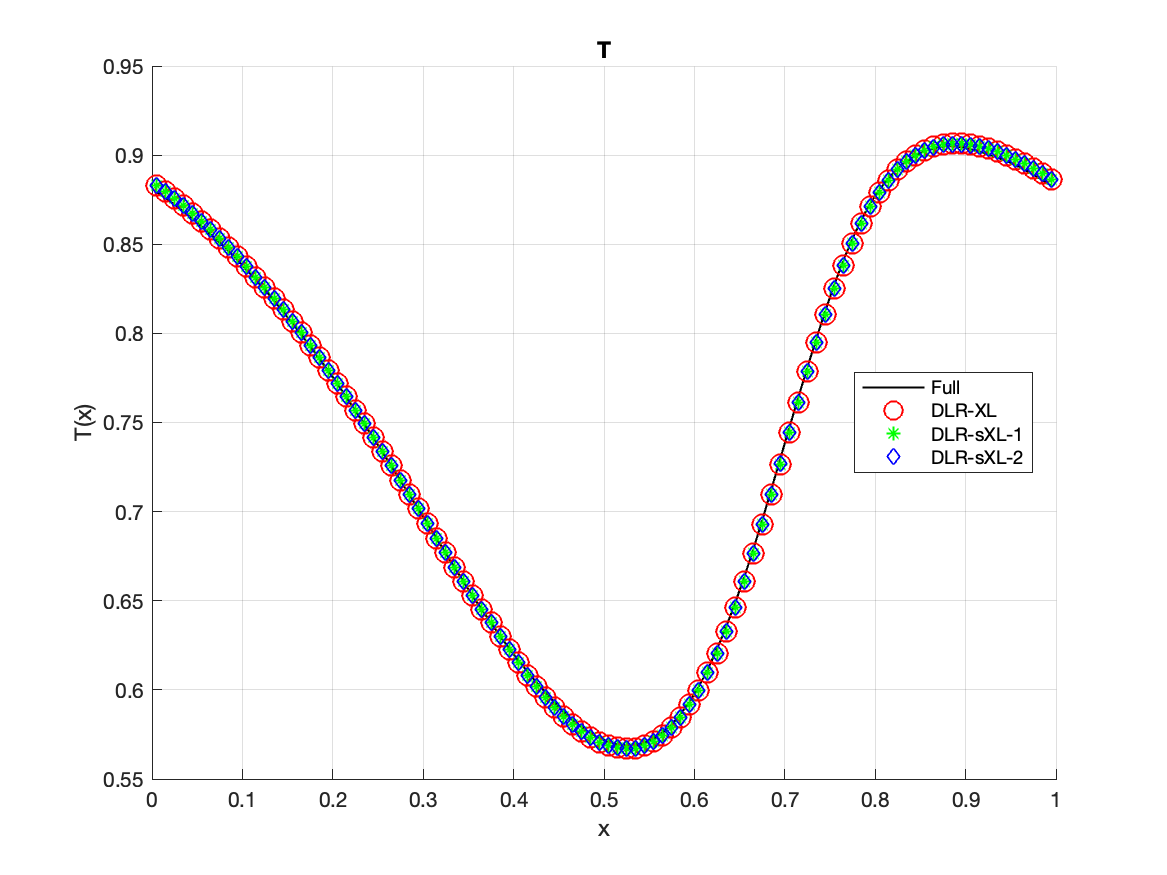}
        \caption{Temperature \(T\).}
    \end{subfigure}
    \begin{subfigure}[b]{0.48\textwidth}
        \includegraphics[width=\textwidth]{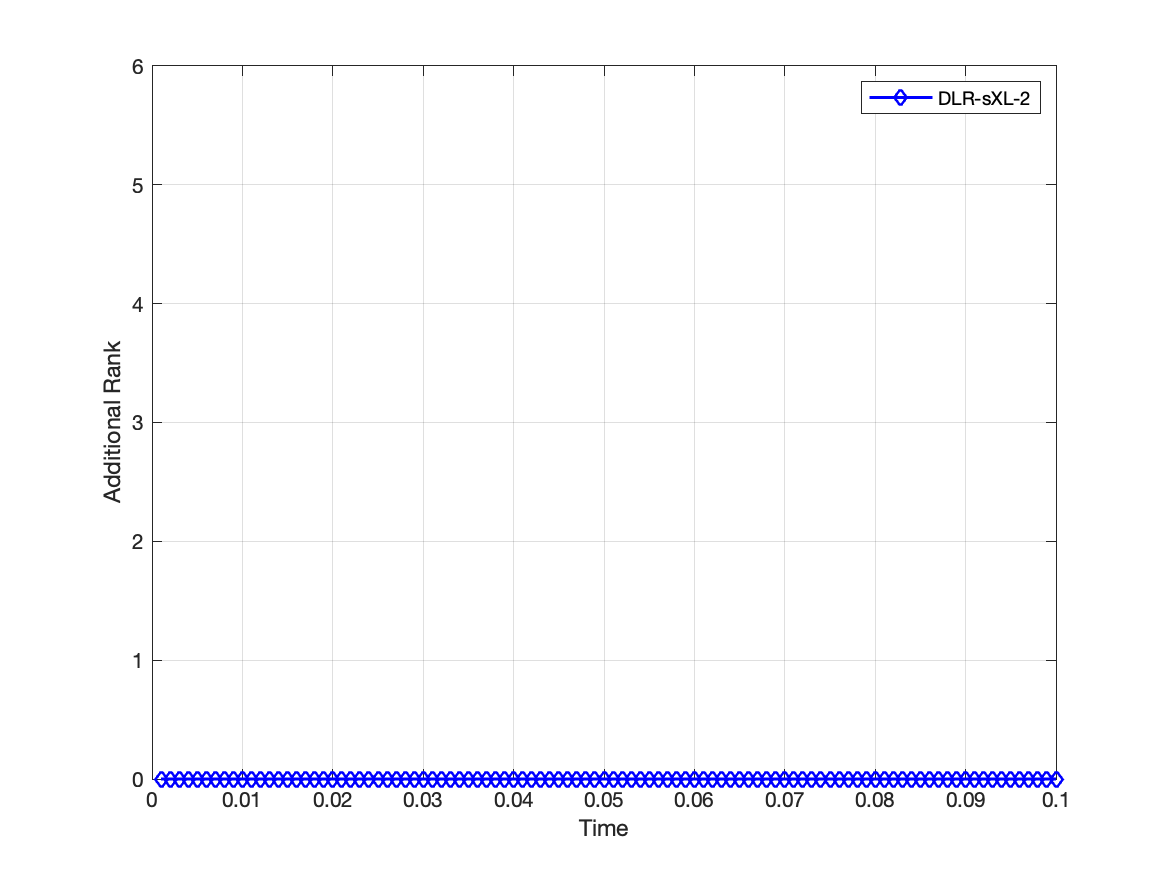}
        \caption{Additional rank of DLR-sXL-2 with \(\text{tol} = 1 \).}
        \label{sin1d}
    \end{subfigure}
    \caption{Comparison of density, bulk velocity (first component), and temperature for the sine problem with \(\Delta t = 10^{-3}\), \(\varepsilon = 1\), $r = 6$ and $\tilde{r} = 6$ for DLR-XL and DLR-sXL-1. The \(\tilde{r}\) of DLR-sXL-2 is shown in \eqref{sin1d} with \(\text{tol} = 1 \). Solid line: full tensor method, red circle: DLR-XL, green star: DLR-sXL-1, blue diamond: DLR-sXL-2.}
    \label{sin1}
\end{figure}

\begin{figure}[htbp]
    \centering
     \begin{subfigure}[b]{0.48\textwidth}
        \includegraphics[width=\textwidth]{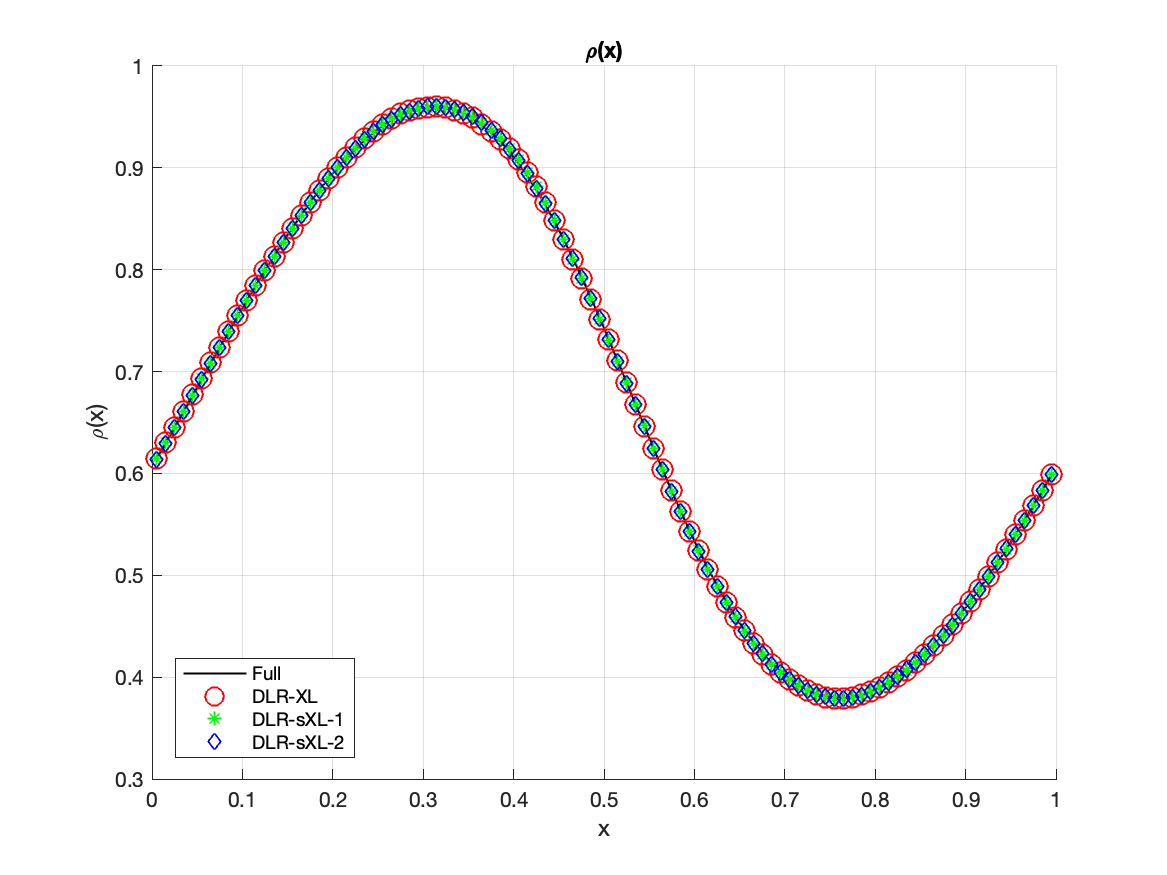}
        \caption{Density \(\rho\).}
    \end{subfigure}
    \begin{subfigure}[b]{0.48\textwidth}
        \includegraphics[width=\textwidth]{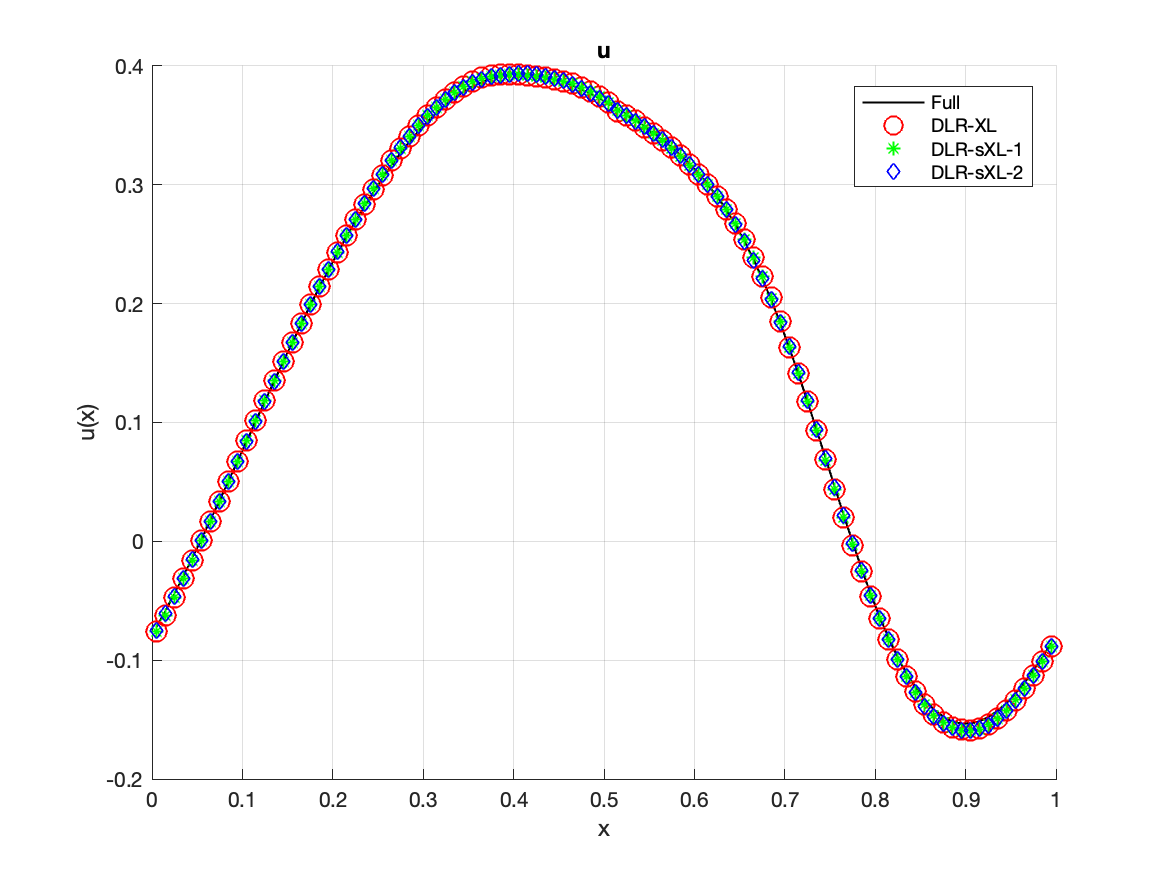}
        \caption{Bulk velocity \(u\).}
    \end{subfigure}
    \begin{subfigure}[b]{0.48\textwidth}
        \includegraphics[width=\textwidth]{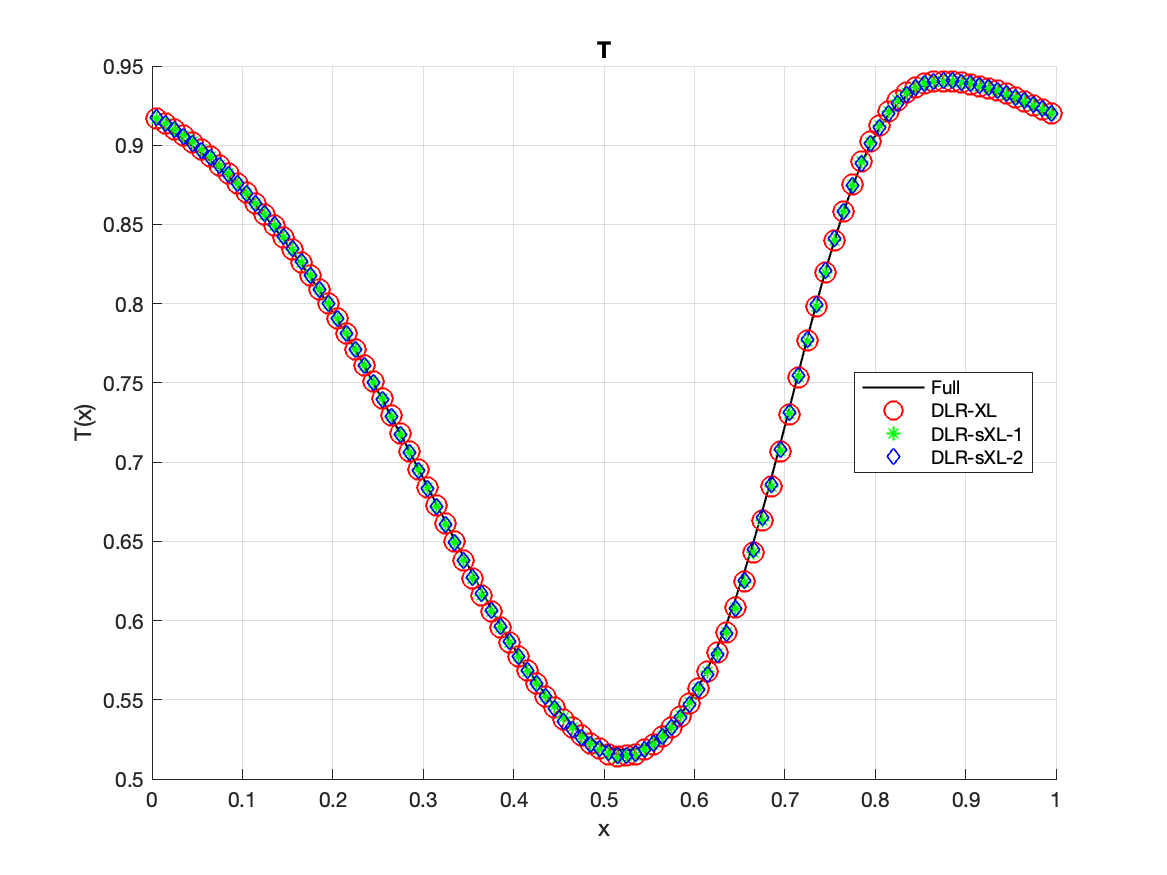}
        \caption{Temperature \(T\).}
    \end{subfigure}
    \begin{subfigure}[b]{0.48\textwidth}
        \includegraphics[width=\textwidth]{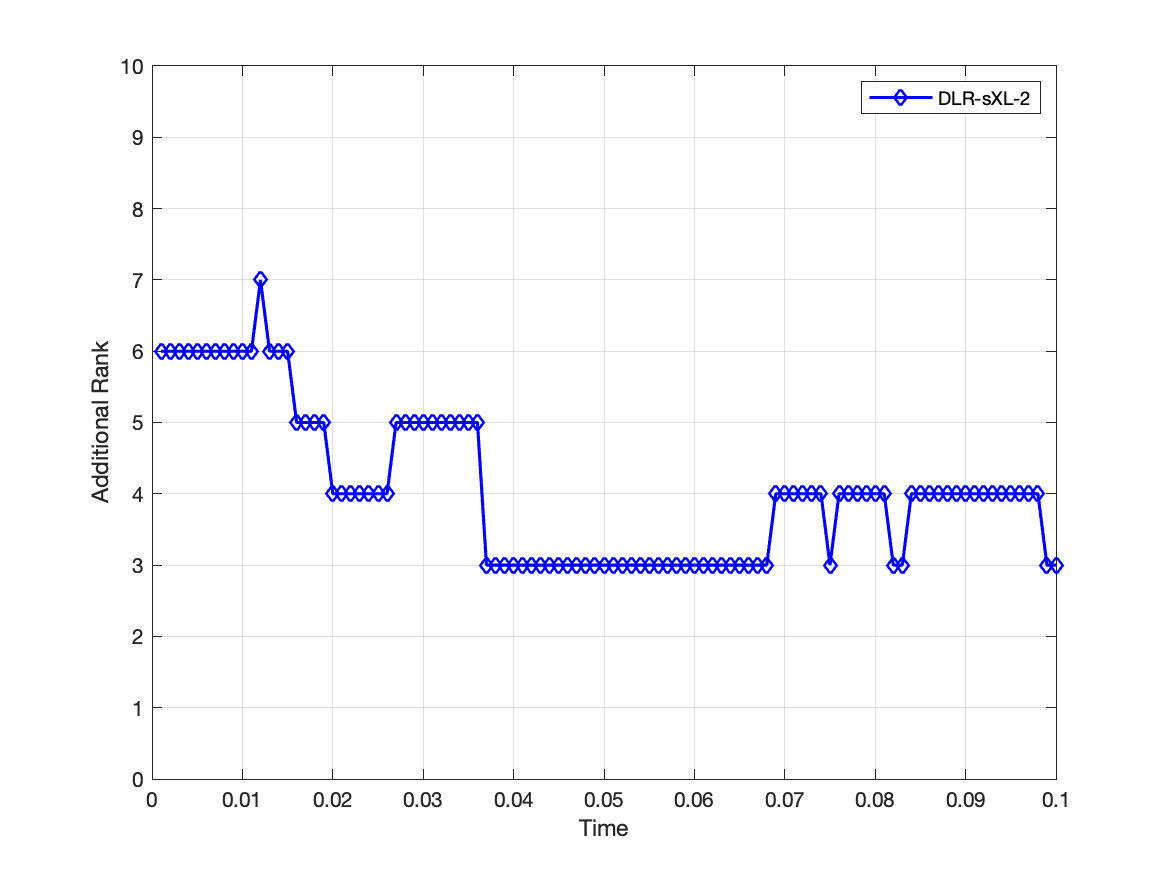}
        \caption{Additional rank of DLR-sXL-2 with \(\text{tol} = 0.05 \).}
        \label{sin1-6d}
    \end{subfigure}
    \caption{Comparison of density, bulk velocity (first component), and temperature for the sine problem with \(\Delta t = 10^{-3}\), \(\varepsilon = 10^{-6}\), $r = 10$ and $\tilde{r} = 10$ for DLR-XL and DLR-sXL-1. The \(\tilde{r}\) of DLR-sXL-2 is shown in \eqref{sin1-6d} with \(\text{tol} = 0.05 \). Solid line: full tensor method, red circle: DLR-XL, green star: DLR-sXL-1, blue diamond: DLR-sXL-2.}
    \label{sin1-6}
\end{figure}

The rank is chosen as \( 6 \) for \(\varepsilon = 1\) and \( 10 \) for \(\varepsilon = 10^{-6}\). The results are presented in Figures \ref{sin1} and \ref{sin1-6}. All low-rank methods—DLR-XL, DLR-sXL-1, and DLR-sXL-2—exhibit comparable accuracy to the full-tensor method. The time step \(\Delta t\) is chosen to be \(1 \times 10^{-3}\), which is substantially larger than \(\varepsilon = 1 \times 10^{-6}\). This choice demonstrates that the time step can be chosen independently of \(\varepsilon\). Additionally, we depict in Figure~\ref{sin1d} and \ref{sin1-6d} the rank augmentation for DLR-sXL-2 with \(\text{tol} = 1\) for $\varepsilon = 1$ and \(\text{tol} = 0.05\) for \(\varepsilon = 10^{-6}\). Notably, when \(\varepsilon = 1\), choosing a sufficiently large tolerance (such that no additional basis functions are added) still yields accurate solutions.

\begin{figure}[htbp]
    \centering
     \begin{subfigure}[b]{0.48\textwidth}
        \includegraphics[width=\textwidth]{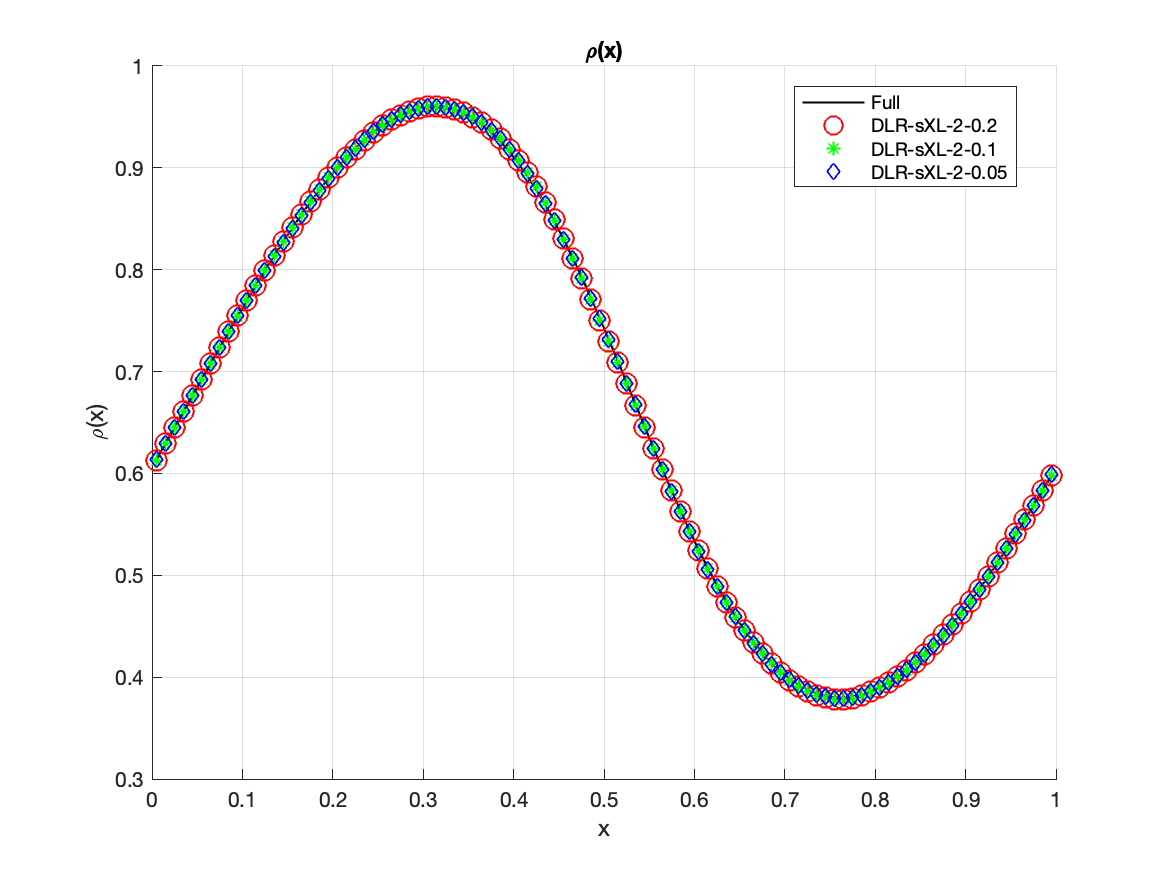}
        \caption{Density \(\rho\).}
    \end{subfigure}
    \begin{subfigure}[b]{0.48\textwidth}
        \includegraphics[width=\textwidth]{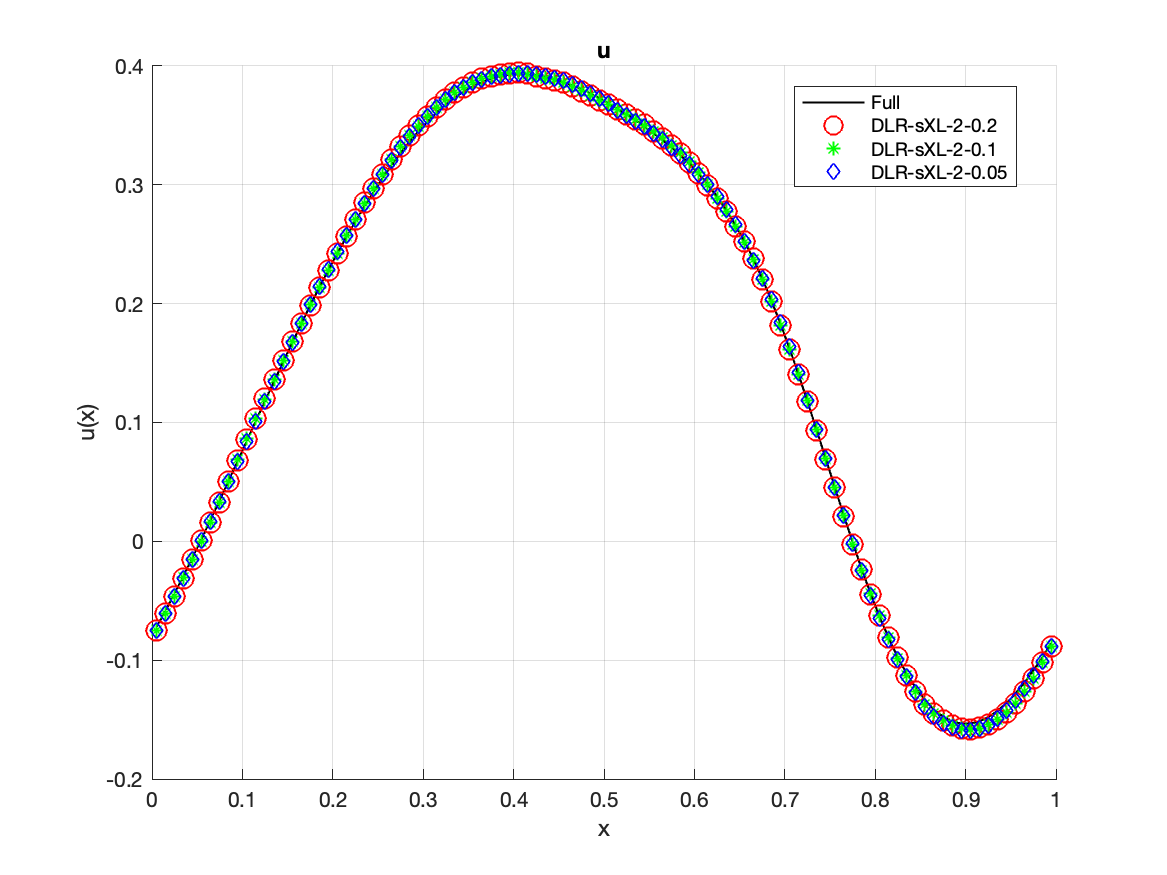}
        \caption{Bulk velocity \(u\).}
    \end{subfigure}
    \begin{subfigure}[b]{0.48\textwidth}
        \includegraphics[width=\textwidth]{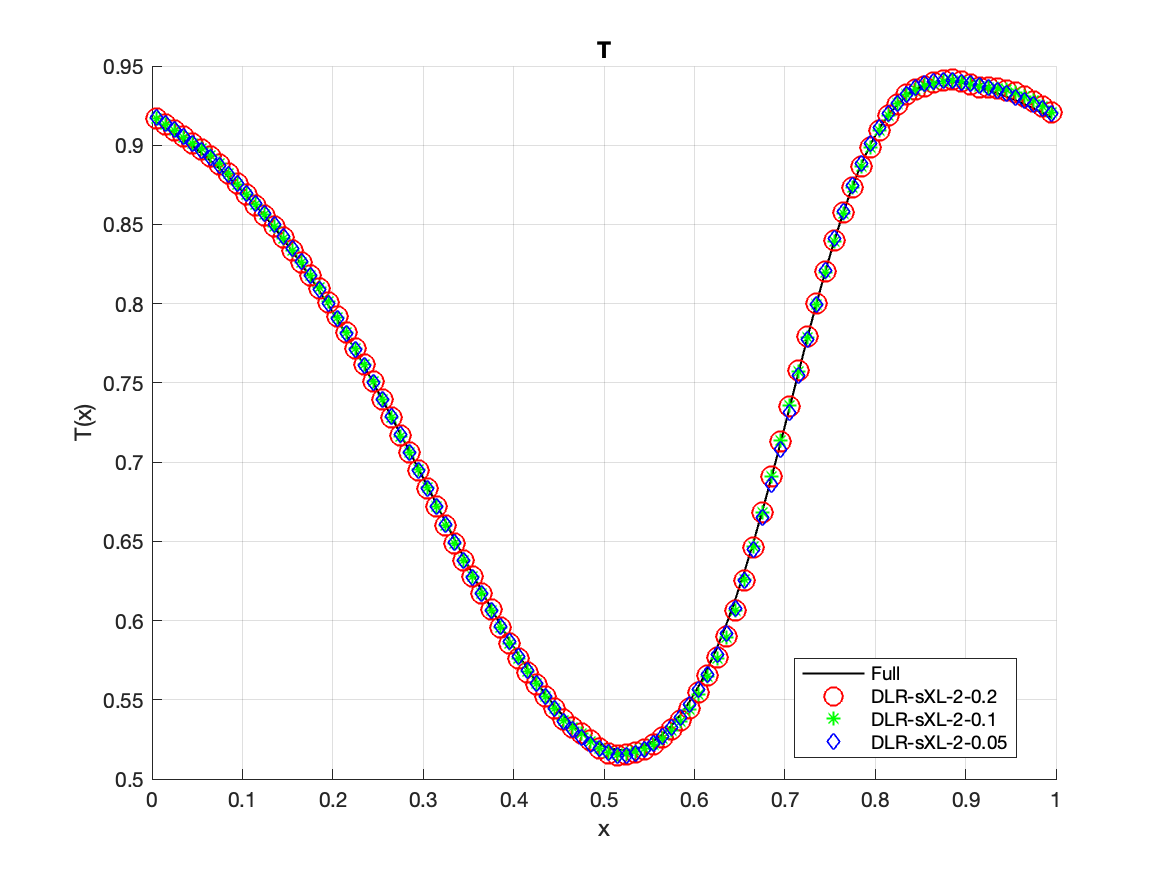}
        \caption{Temperature \(T\).}
        \label{sin1-6-3c}
    \end{subfigure}
    \begin{subfigure}[b]{0.48\textwidth}
        \includegraphics[width=\textwidth]{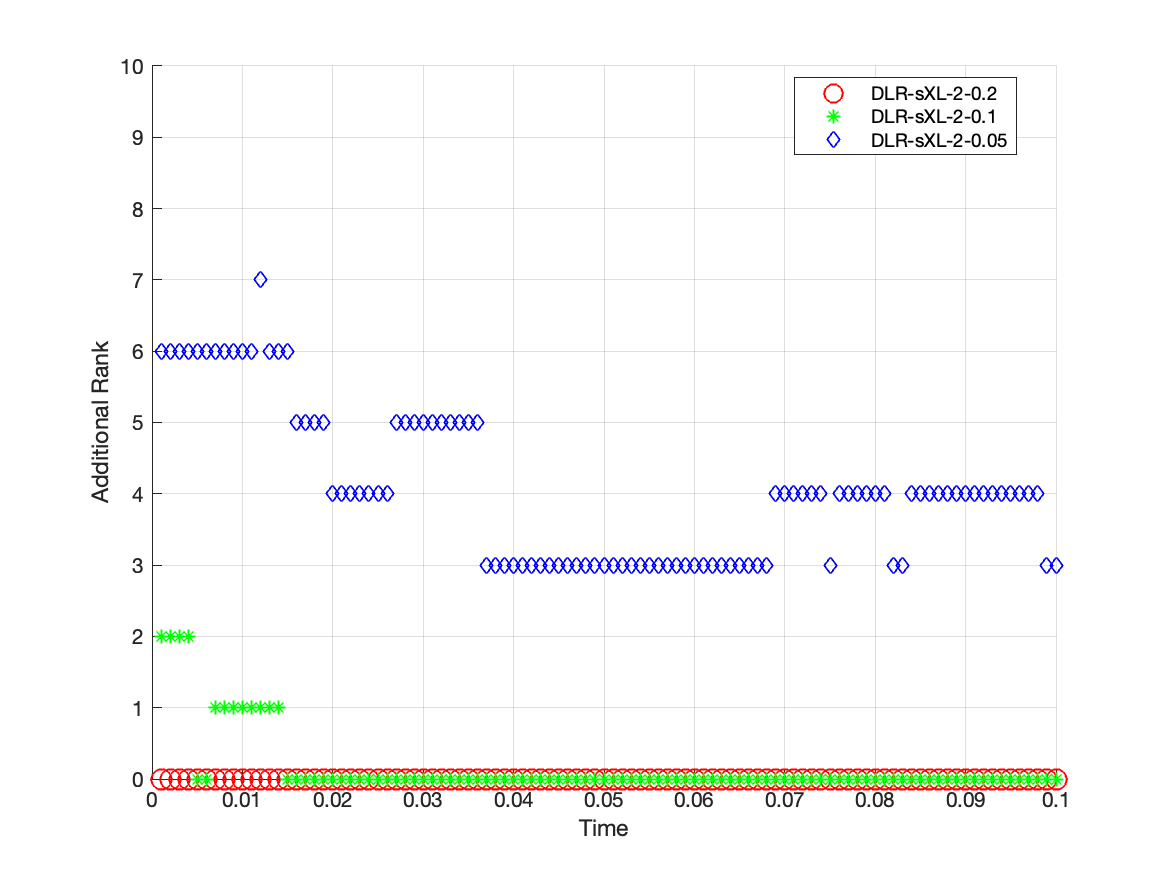}
        \caption{Additional ranks of DLR-sXL-2 with different $\text{tols}$.}
        \label{sin1-6-3d}
    \end{subfigure}
    \caption{Comparison of DLR-sXL-2 with rank \(r = 10\) under three different tolerances, \(\text{tol} \in \{0.2, 0.1, 0.05\}\), for \(\varepsilon = 10^{-6}\). Solid line: full tensor method, red circle: DLR-sXL-2 with $\text{tol} = 0.2$, green star: DLR-sXL-2 with $\text{tol} = 0.1$, blue diamond: DLR-sXL-2 with $\text{tol} = 0.05$. In Figure~\ref{sin1-6-3c}, one can observe a slight discrepancy between DLR-sXL-2-0.2 and others. The additional ranks refers to the number of extra basis functions introduced during the augmentation steps.}
    \label{sin1-6-3}
\end{figure}

To investigate the effect of different tolerances in DLR-sXL-2 for small $\varepsilon = 10^{-6}$, we present Figure~\ref{sin1-6-3}, where three values of \(\text{tol}\), \(0.2\), \(0.1\), and \(0.05\), are tested. For \(\text{tol} = 0.2\), no additional basis functions is added during augmentation step, yet the solution changes only slightly (see Figure~\ref{sin1-6-3c}). When \(\text{tol} = 0.1\), only a few basis functions are added at the first few steps, and the overall accuracy remains high. Finally, for \(\text{tol} = 0.05\), more basis functions are added, but the solution accuracy does not improve noticeably. These observations imply that the augmentation step is unnecessary for large \(\varepsilon\), while for small \(\varepsilon\), adding fewer than \(\tilde{r}\) basis functions can still produce accurate solutions.


\subsection{Shock tube problem}
We then examine the performance of the proposed integrators for discontinuous initial data using the classical shock tube problem. 

The initial data is at the equilibrium: 
\begin{equation}
    f^0(x, v, z) = \frac{\rho^0}{2\pi T^0} e^{-\frac{|v-u^0|^2}{2T^0}},
\end{equation}
where the density, bulk velocity and temperature are given by
\[
\rho^0(x) = 
\begin{cases}
1, & x \le 0.5, \\
0.125, & x > 0.5,
\end{cases}
\quad
u^0(x) = 
\begin{cases}
(0,0), & x \le 0.5, \\
(0,0), & x > 0.5,
\end{cases}
\quad
T^0(x) =
\begin{cases}
1, & x \le 0.5, \\
0.25, & x > 0.5.
\end{cases}
\]

\begin{figure}[htbp]
    \centering
     \begin{subfigure}[b]{0.48\textwidth}
        \includegraphics[width=\textwidth]{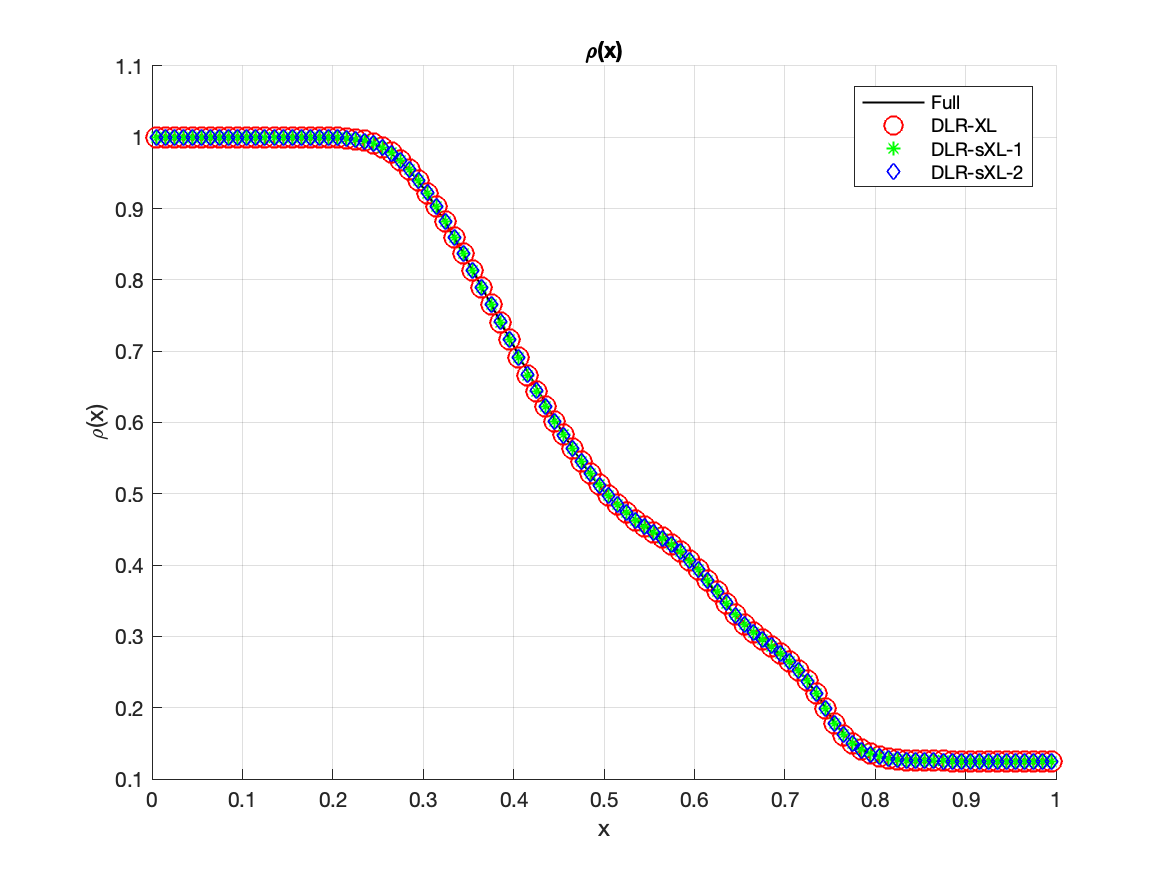}
        \caption{Density \(\rho\).}
    \end{subfigure}
    \begin{subfigure}[b]{0.48\textwidth}
        \includegraphics[width=\textwidth]{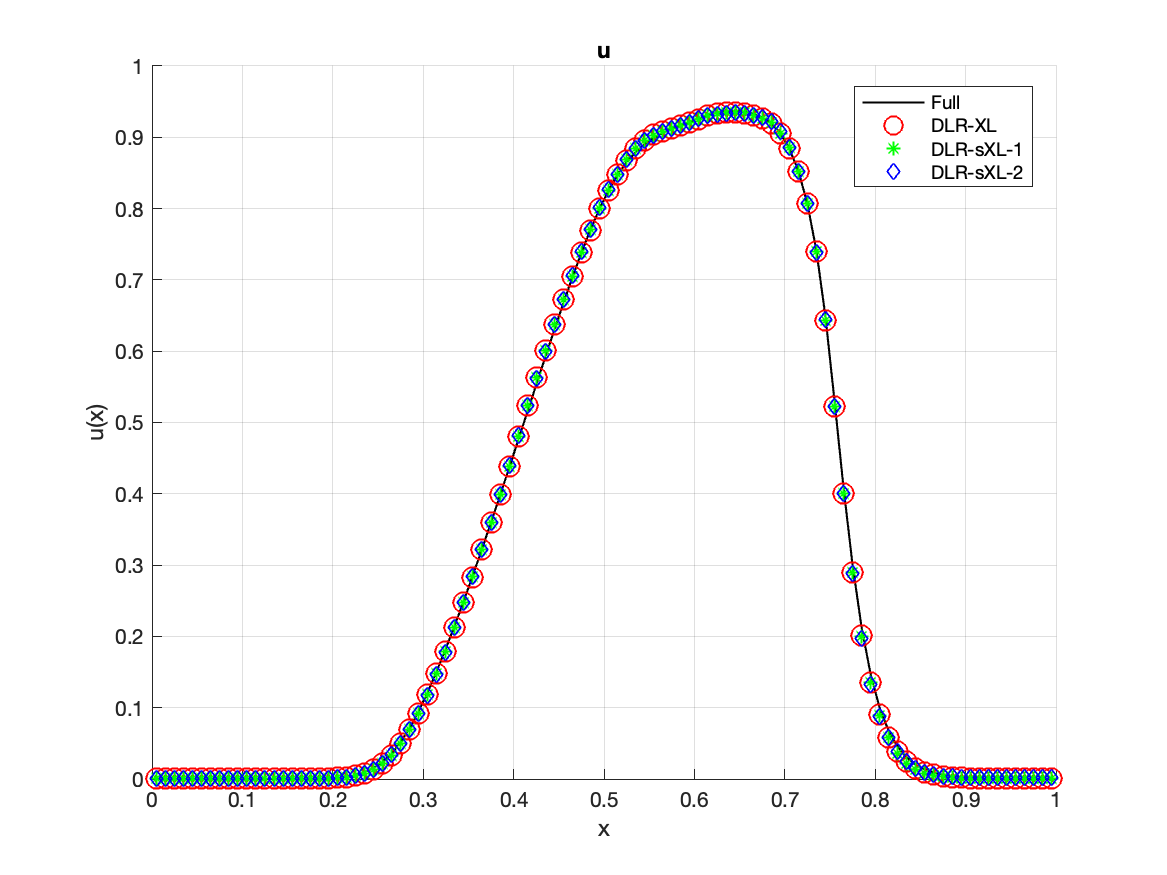}
        \caption{Bulk velocity \(u\).}
    \end{subfigure}
    \begin{subfigure}[b]{0.48\textwidth}
        \includegraphics[width=\textwidth]{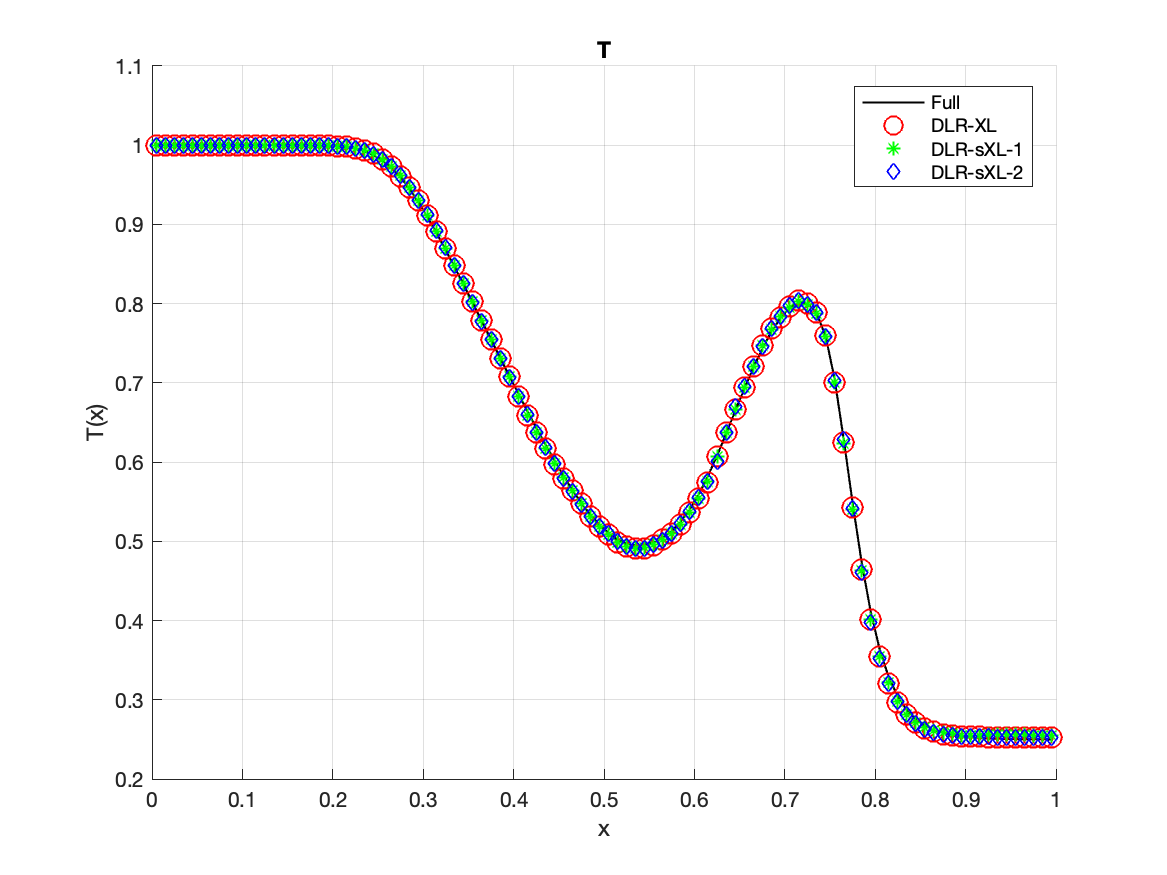}
        \caption{Temperature \(T\).}
    \end{subfigure}
    \begin{subfigure}[b]{0.48\textwidth}
        \includegraphics[width=\textwidth]{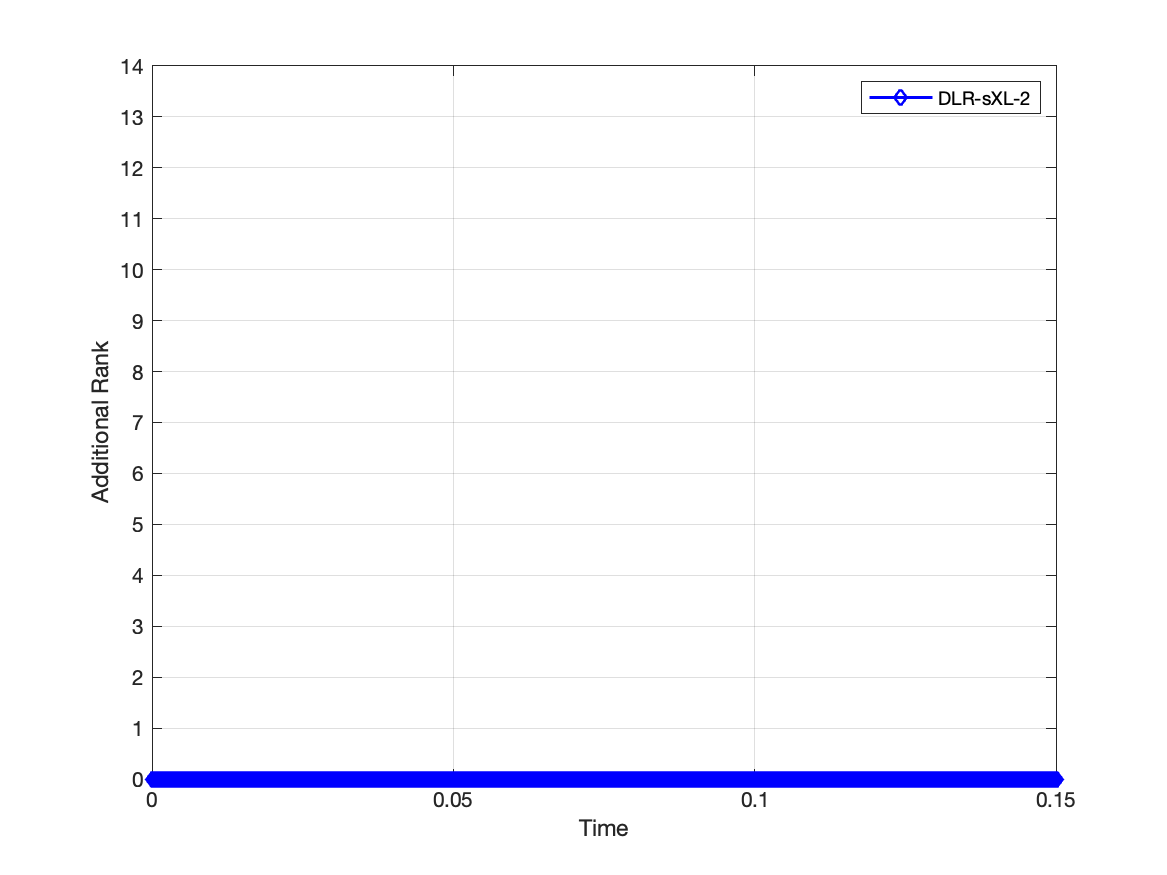}
        \caption{Additional rank of DLR-sXL-2 with \(\text{tol} = 1 \).}
        \label{shock1-3d}
    \end{subfigure}
    \caption{Comparison of density, bulk velocity (first component), and temperature for the shock tube problem with \(\Delta t = 10^{-4}\), \(\varepsilon = 10^{-3}\), $r = 14$ and $\tilde{r} = 14$ for DLR-XL and DLR-sXL-1. The \(\tilde{r}\) of DLR-sXL-2 is shown in \eqref{shock1-3d} with \(\text{tol} = 1 \). Solid line: full tensor method, red circle: DLR-XL, green star: DLR-sXL-1, blue diamond: DLR-sXL-2.}
    \label{shock1-3}
\end{figure}

\begin{figure}[htbp]
    \centering  
     \begin{subfigure}[b]{0.48\textwidth}
        \includegraphics[width=\textwidth]{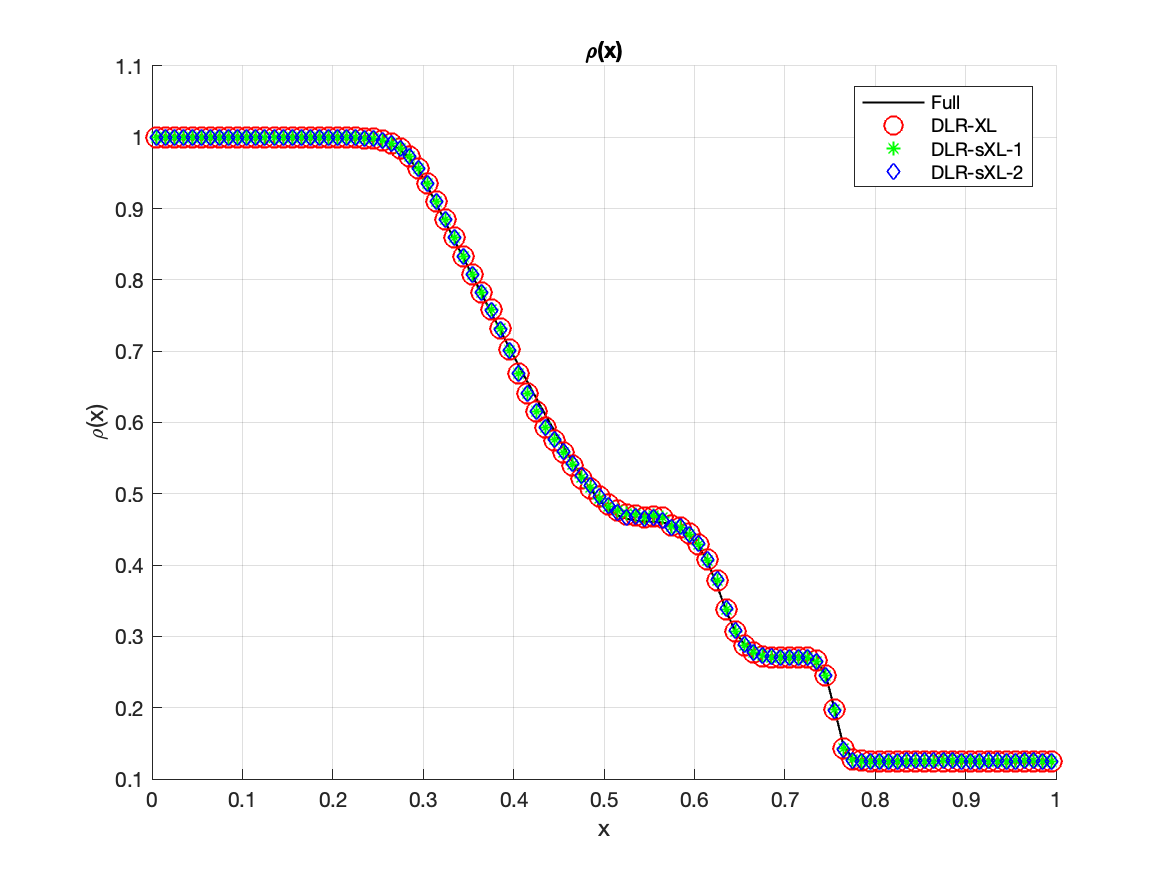}
        \caption{Density \(\rho\).}
    \end{subfigure}
    \begin{subfigure}[b]{0.48\textwidth}
        \includegraphics[width=\textwidth]{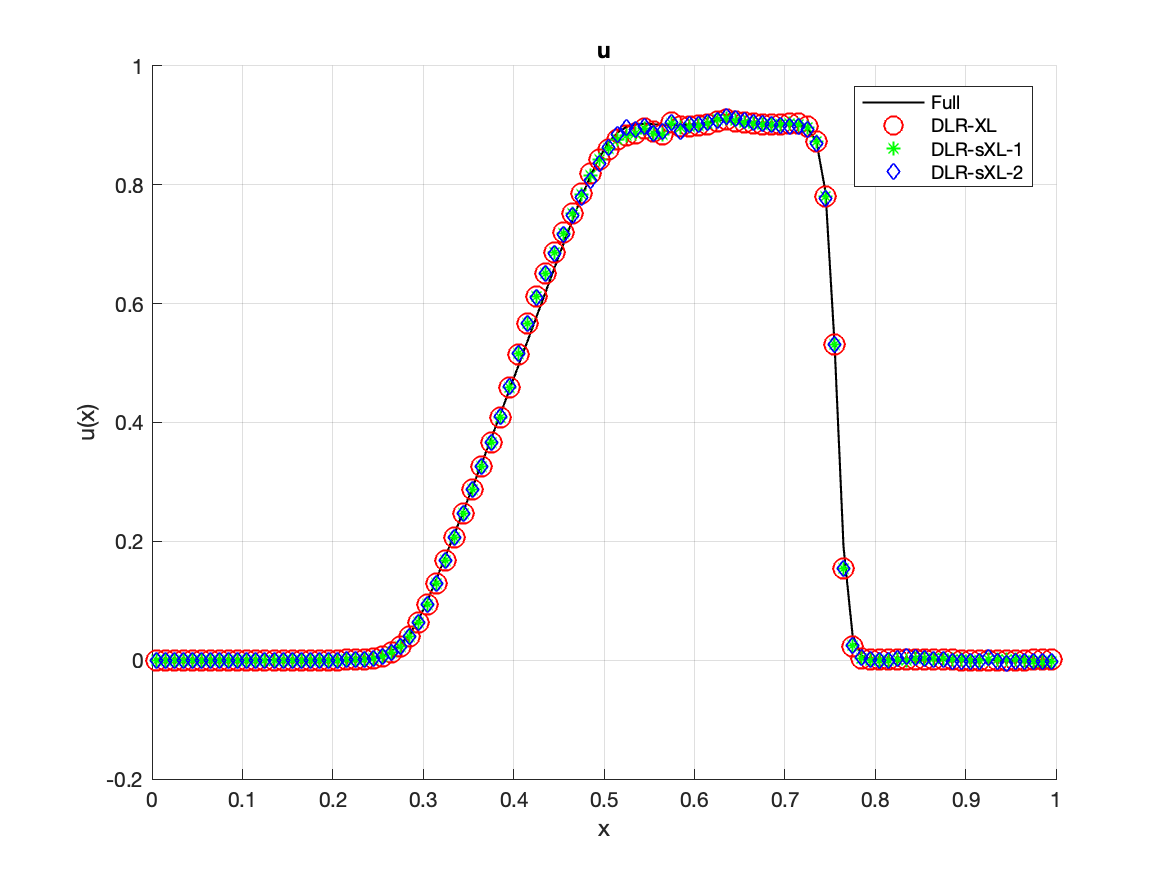}
        \caption{Bulk velocity \(u\).}
    \end{subfigure}
    \begin{subfigure}[b]{0.48\textwidth}
        \includegraphics[width=\textwidth]{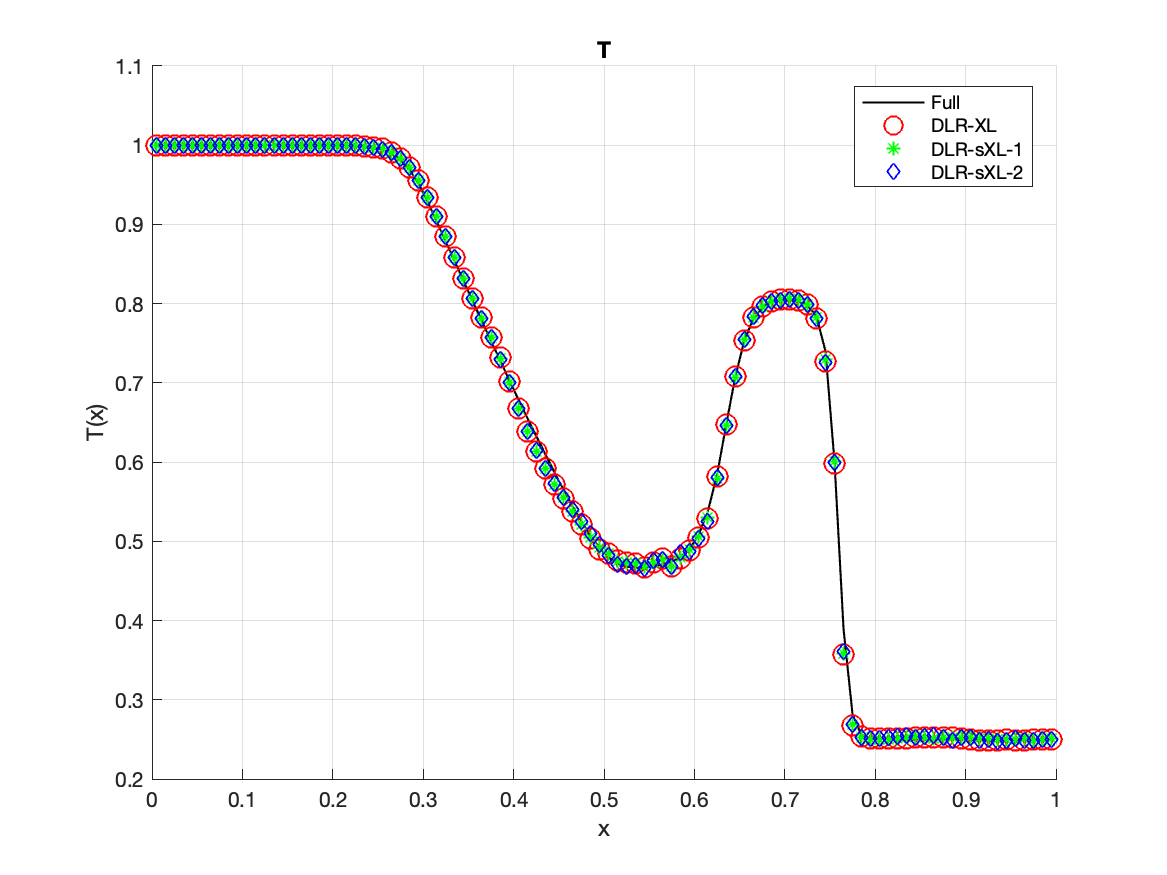}
        \caption{Temperature \(T\).}
    \end{subfigure}
    \begin{subfigure}[b]{0.48\textwidth}
        \includegraphics[width=\textwidth]{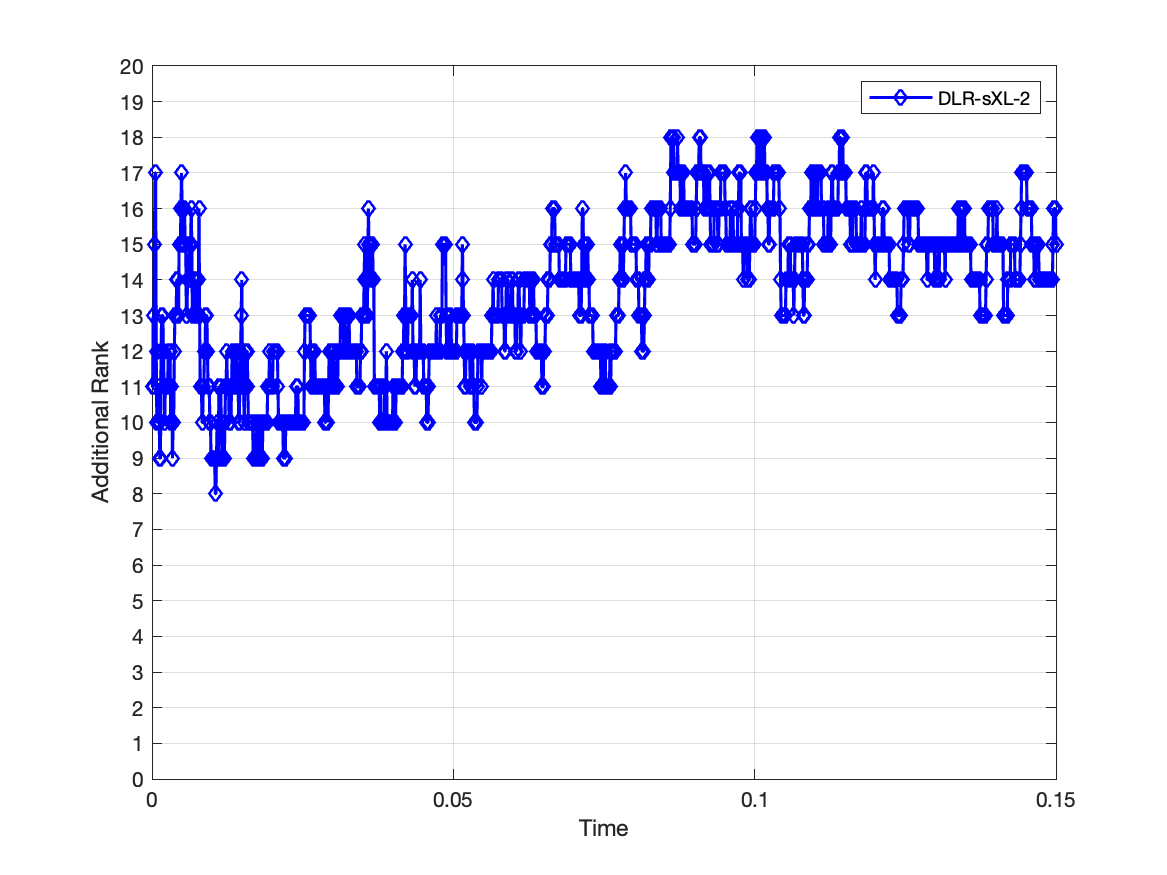}
        \caption{Additional rank of DLR-sXL-2 with \(\text{tol} = 0.01 \).}
        \label{shock1-6d}
    \end{subfigure}
    \caption{Comparison of density, bulk velocity (first component), and temperature for the shock tube problem with \(\Delta t = 10^{-4}\), \(\varepsilon = 10^{-6}\), $r = 20$ and $\tilde{r} = 20$ for DLR-XL and DLR-sXL-1. The \(\tilde{r}\) of DLR-sXL-2 is shown in \eqref{shock1-6d} with \(\text{tol} = 0.01 \). Solid line: full tensor method, red circle: DLR-XL, green star: DLR-sXL-1, blue diamond: DLR-sXL-2.}
    \label{shock1-6}
\end{figure}

We set the rank to \(14\) for \(\varepsilon = 10^{-3}\) and to \(20\) for \(\varepsilon = 10^{-6}\); see Figures~\ref{shock1-3} and \ref{shock1-6}. To compare with other methods, the tolerance \(\text{tol}\) is chosen as \(1\) when \(\varepsilon = 10^{-3}\) and \(0.01\) when \(\varepsilon = 10^{-6}\). In both cases, all low-rank methods exhibit accuracy comparable to the full-tensor method. The time step \(\Delta t = 1 \times 10^{-4}\) is selected solely by the CFL condition, independent of \(\varepsilon\).

Figures~\ref{shock1-3d} and \ref{shock1-6d} illustrate the rank augmentation for DLR-sXL-2 using \(\text{tol} = 1\) and \(\text{tol} = 0.01\) under different values of \(\varepsilon\). Notably, when \(\varepsilon = 1\) and \(\text{tol}\) is large enough that no additional basis functions are introduced, DLR-sXL-2 still achieves accurate solutions. 

\begin{figure}[htbp]
    \centering  
     \begin{subfigure}[b]{0.48\textwidth}
        \includegraphics[width=\textwidth]{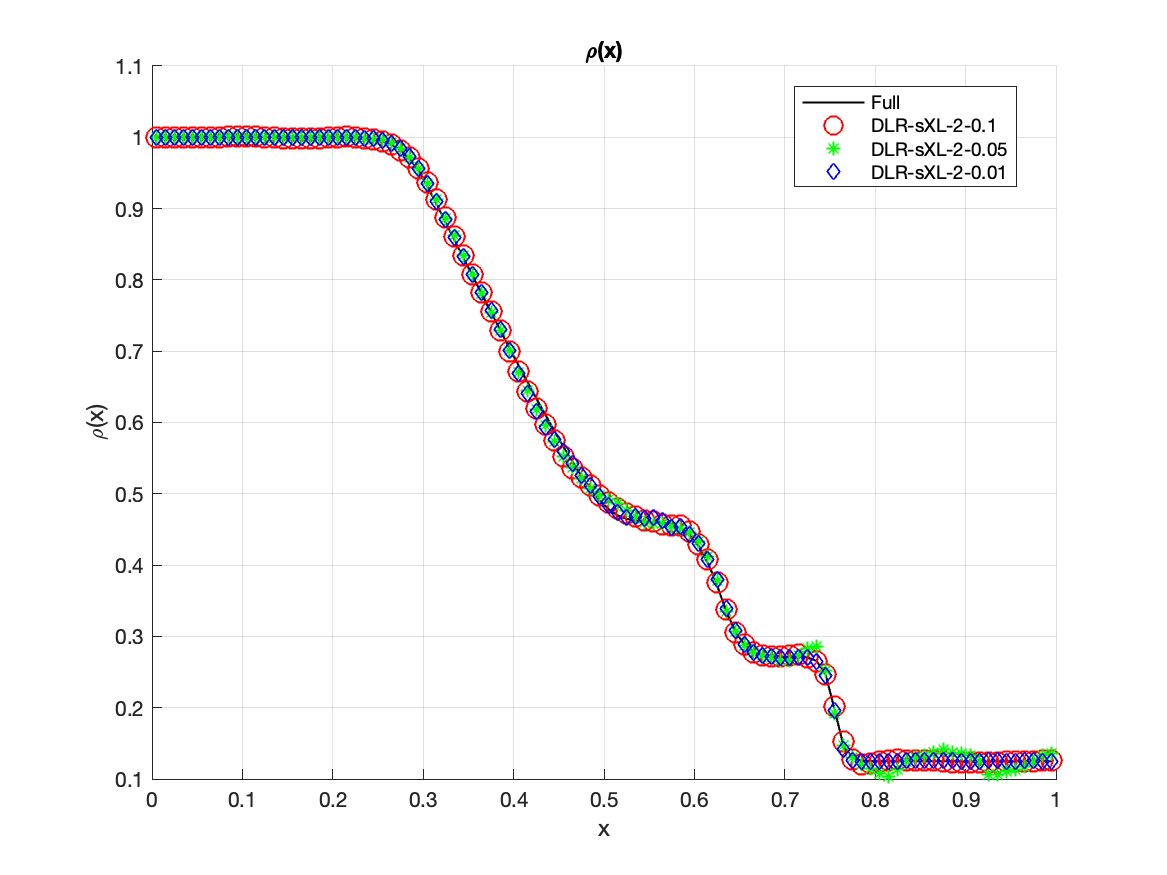}
        \caption{Density \(\rho\).}
        \label{shock1-6-3a}
    \end{subfigure}
    \begin{subfigure}[b]{0.48\textwidth}
        \includegraphics[width=\textwidth]{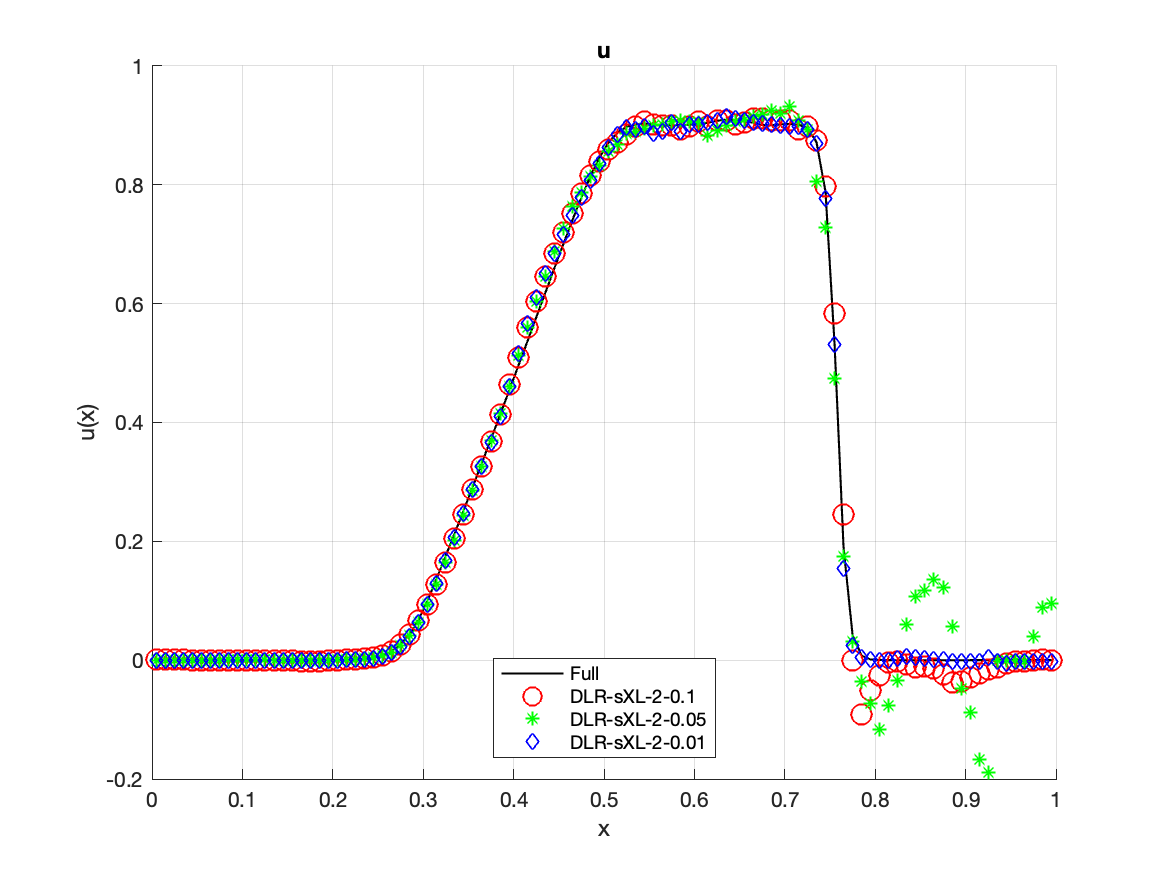}
        \caption{Bulk velocity \(u\).}
        \label{shock1-6-3b}
    \end{subfigure}
    \begin{subfigure}[b]{0.48\textwidth}
        \includegraphics[width=\textwidth]{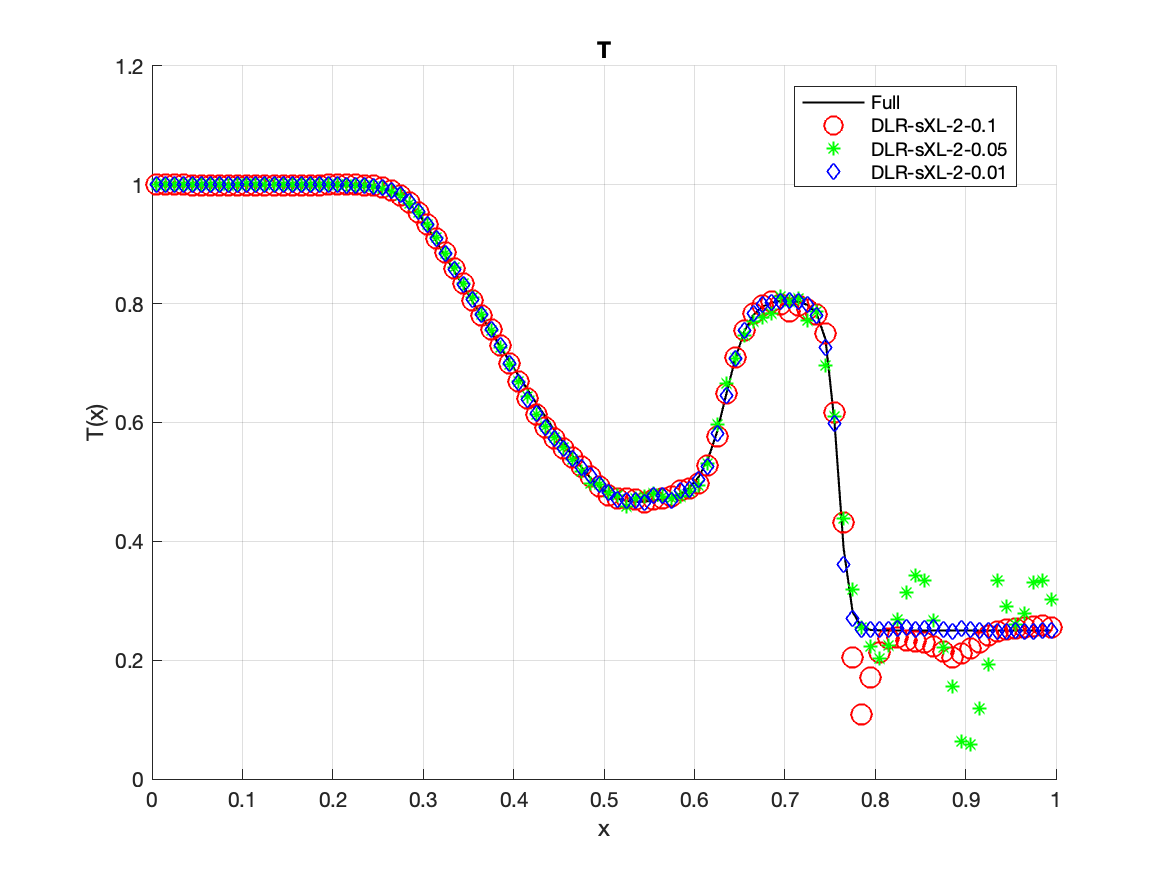}
        \caption{Temperature \(T\).}
        \label{shock1-6-3c}
    \end{subfigure}
    \begin{subfigure}[b]{0.48\textwidth}
        \includegraphics[width=\textwidth]{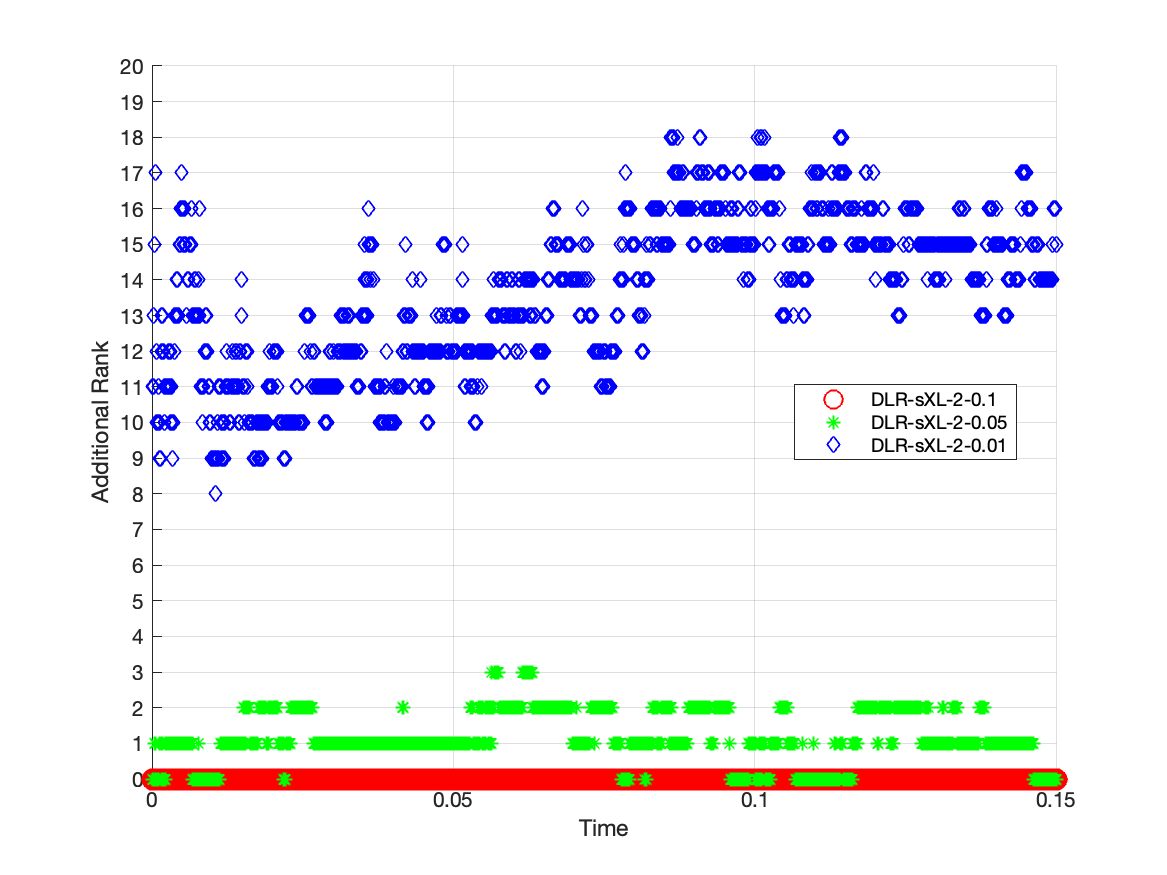}
        \caption{Additional ranks of DLR-sXL-2 with different $\text{tols}$.}
        \label{shock1-6-3d}
    \end{subfigure}
    \caption{Comparison of DLR-sXL-2 with rank \(r = 20\) under three different tolerances, \(\text{tol} \in \{0.1, 0.05, 0.01\}\), for \(\varepsilon = 10^{-6}\). Solid line: full tensor method, red circle: DLR-sXL-2 with $\text{tol} = 0.1$, green star: DLR-sXL-2 with $\text{tol} = 0.05$, blue diamond: DLR-sXL-2 with $\text{tol} = 0.01$. The additional ranks refers to the number of extra basis functions introduced during the augmentation steps. DLR-sXL-2 performs poorly with larger tolerances $\text{tol} = 0.1$ and 0.05.}
    \label{shock1-6-3}
\end{figure}

To examine the effect of varying \(\text{tol}\) in DLR-sXL-2 for the small \(\varepsilon = 10^{-6}\), we refer to Figure~\ref{shock1-6-3}, where three \text{tols}, \(0.1\), \(0.05\), and \(0.01\), are tested. For \(\text{tol} = 0.1\) and \(0.05\), few or no basis functions are added, and the method performs poorly (see Figures~\ref{shock1-6-3a}, \ref{shock1-6-3b} and \ref{shock1-6-3c}). By contrast, \(\text{tol} = 0.01\) results in more basis functions being introduced and maintains high accuracy. Even so, the augmented rank remains below \(\tilde{r} = 20\). These observations suggest that, when DLR-XL or DLR-sXL-1 performs reliably, one can choose an appropriate \(\text{tol}\) that adds fewer than \(\tilde{r}\) basis functions while still ensuring accurate solutions.


\section{Conclusions}
\label{sec:conclusions}

We have demonstrated that an efficient dynamical low-rank algorithm for the stiff Boltzmann equation can be obtained using the proposed (s)XL integrator. In particular, the algorithm requires evaluating the Boltzmann collision operator only $r^2$ times per time step, where $r$, the rank of the approximation, is much smaller than the number of spatial grid points. Moreover, the scheme is asymptotic-preserving without requiring either implicit treatment of the Boltzmann collision operator or expensive evaluation of the Maxwellian. We have conducted a number of 1+2 dimensional simulations that confirm the favorable performance of the method. 
We 
note that the XL and sXL integrators considered here are first-order accurate in time. This is in fact limited by the first-order accuracy of the full tensor AP scheme we are based on. Extension to higher accuracy in time presents interesting future work.


\section*{Acknowledgments}

JH and SZ's research is partially supported by NSF DMS-2409858, DoD MURI FA9550-24-1-0254, and DOE DE-SC0023164.

\bibliographystyle{siamplain}
\bibliography{hu_bibtex,einkemmer_bibtex}
\end{document}